\documentclass[11pt]{article}

\usepackage{graphicx}

\usepackage{labelfig}

\usepackage{amsmath,amsfonts,amssymb,amsthm,upref,bm}%

\usepackage{mathtools}

\usepackage{mathrsfs}

\usepackage{booktabs}

\usepackage[paperwidth=210mm,paperheight=297mm,top=20mm,bottom=20mm,left=20mm,right=20mm]{geometry}

\DeclareMathAlphabet{\varmathbb}{U}{pxsyb}{m}{n}
\DeclareMathAlphabet{\mathpzc}{OT1}{pzc}{m}{it}

\newcommand{\MF}[1]{\mathop{#1}\nolimits}

\newcommand{\ii}{\kern0.05em\mathrm{i}\kern0.05em} \newcommand{\E}[1]{\textrm{e}^{#1}}

\renewcommand{\vec}[1]{\bm{#1}}

\renewcommand{\Re}{\MF{\mathrm{Re}}} \renewcommand{\Im}{\MF{\mathrm{Im}}}

\newcommand{\sign}{\mathop{\textrm{sign}}}

\newcommand{\Zzo}{\varmathbb{Z}_{\scriptscriptstyle0}^{-}}%
\newcommand{\Cb}{\varmathbb{C}}%

\newcommand{\Pcoef}{\mathcal{P}}
\newcommand{\Qcoef}{\mathcal{Q}}
\newcommand{\Rcoef}{\mathcal{R}}
\newcommand{\Scoef}{\mathcal{S}}

\newcommand{\SectInf}{S}%

\newcommand{\BCzero}{\mathscr{B}_{0\infty}}%
\newcommand{\BCone}{\mathscr{B}_{1\infty}}%

\def\midpoint#1{z^\star_{#1}}%
\def\nsp{n_*}%
\def\Lkr{\mathscr{L}}%
\def\reservecoeff{\varkappa}
\def\car{\mathpzc{c}}

\def\HeunC{\displaystyle\setbox0=\hbox{\tiny\itshape\bfseries \kern1.52mm c}\wd0=0pt\ht0=0pt\dp0=0pt\MF{\textit{\hbox{\raise1.67ex\copy0}H\/}}}%

\def\HeunCtxst{\scalebox{0.735}{\ensuremath{\HeunC}}}%

\def\HeunCl{\setbox0=\hbox{\tiny\itshape\bfseries \kern1.52mm c}\wd0=0pt\ht0=0pt\dp0=0pt\MF{\textit{\hbox{\raise1.67ex\copy0}Hl\/}}}%
\def\HeunCltxst{\scalebox{0.735}{\ensuremath{\HeunCl}}}%

\def\HeunCs{\setbox0=\hbox{\tiny\itshape\bfseries \kern1.52mm c}\wd0=0pt\ht0=0pt\dp0=0pt\MF{\textit{\hbox{\raise1.67ex\copy0}Hs\/}}}%

\def\psHeunC{\setbox0=\hbox{\tiny\itshape\bfseries \kern1.56mm c}\wd0=0pt\ht0=0pt\dp0=0pt\MF{\hbox{\raise1.67ex\copy0}\mathpzc{H}}}%

\def\psHeunCo{\setbox0=\hbox{\tiny\itshape\bfseries \kern1.56mm c}\wd0=0pt\ht0=0pt\dp0=0pt%
\setbox1=\hbox{\kern3.475mm\scalebox{0.5}{$\bm{\circ}$}}\wd1=0pt\ht1=0pt\dp1=0pt%
\MF{\hbox{\raise1.67ex\copy0}\hbox{\raise0.95ex\copy1}\mathpzc{H}\kern2pt}}%

\def\psHeunCl{\setbox0=\hbox{\tiny\itshape\bfseries \kern1.56mm c}\wd0=0pt\ht0=0pt\dp0=0pt\MF{\hbox{\raise1.67ex\copy0}\mathpzc{Hl}_{\!0}}}%

\def\psHeunClpath#1{\setbox0=\hbox{\tiny\itshape\bfseries \kern1.52mm c}\wd0=0pt\ht0=0pt\dp0=0pt\MF{\hbox{\raise1.67ex\copy0}\mathpzc{Hl}_{#1}}}%
\def\psHeunCspath#1{\setbox0=\hbox{\tiny\itshape\bfseries \kern1.52mm c}\wd0=0pt\ht0=0pt\dp0=0pt\MF{\hbox{\raise1.67ex\copy0}\mathpzc{Hs}_{#1}}}%

\def\psHeunClfin{\setbox0=\hbox{\tiny\itshape\bfseries \kern1.52mm c}\wd0=0pt\ht0=0pt\dp0=0pt\MF{\hbox{\raise1.67ex\copy0}\mathpzc{Hl}}}%
\def\psHeunCsfin{\setbox0=\hbox{\tiny\itshape\bfseries \kern1.52mm c}\wd0=0pt\ht0=0pt\dp0=0pt\MF{\hbox{\raise1.67ex\copy0}\mathpzc{Hs}}}%

\def\psHeunClnearone{\setbox0=\hbox{\tiny\itshape\bfseries \kern1.52mm c}\wd0=0pt\ht0=0pt\dp0=0pt\MF{\hbox{\raise1.67ex\copy0}\mathpzc{Hl}^{(1)}}}%
\def\psHeunCsnearone{\setbox0=\hbox{\tiny\itshape\bfseries \kern1.52mm c}\wd0=0pt\ht0=0pt\dp0=0pt\MF{\hbox{\raise1.67ex\copy0}\mathpzc{Hs}^{(1)}}}%

\def\psHeunClnearinf{\setbox0=\hbox{\tiny\itshape\bfseries \kern1.52mm c}\wd0=0pt\ht0=0pt\dp0=0pt\MF{\hbox{\raise1.67ex\copy0}\mathpzc{Hl}^{(\infty)}}}%
\def\psHeunCsnearinf{\setbox0=\hbox{\tiny\itshape\bfseries \kern1.52mm c}\wd0=0pt\ht0=0pt\dp0=0pt\MF{\hbox{\raise1.67ex\copy0}\mathpzc{Hs}^{(\infty)}}}%

\allowdisplaybreaks[3]

\begin{document}

\title{On evaluation of the confluent Heun functions}

\author{Oleg\ V.\ Motygin}
 
\date{Institute of Problems in Mechanical Engineering, Russian Academy\\ of
Sciences,
 V.O., Bol'shoy pr., 61, 199178 St.\,Petersburg, Russia\\ 
email: \texttt{o.v.motygin@gmail.com}}

\maketitle

\begin{abstract}

In this paper we consider the confluent Heun equation, which is a linear differential equation of second order with three singular points --- two of them are regular and the third one is irregular of rank 1. The purpose of the work is to propose a procedure for numerical evaluation of the equation's solutions (confluent Heun functions). A scheme based on power series, asymptotic expansions and analytic continuation is described. Results of numerical tests are given.

\end{abstract}

\section{Introduction}
\label{intro0}

Heun differential equation was introduced by Karl Heun in 1889 \cite{Heun} as a generalization of the hypergeometric one. The general Heun equation is a Fuchsian equation with four regular singular points, which are usually chosen to be $z = 0$, $1$, $a$, and $\infty$ in the complex $z$-plane.
Various kinds of confluence of these singularities, when two or more of them merge to an irregular singularity, produce the confluent, double confluent, biconfluent and triconfluent Heun equations. For a comprehensive mathematical treatment of the topic, we refer to \cite{Ronveaux1995,SlavyanovLay2000,SleemanKuznetzov2010}. In the present paper, we deal with the confluent Heun equation being a result of the simplest case of confluence $a\to\infty$ and having two regular singular points $z=0$, $1$ and an irregular one $z=\infty$.  

The solutions of the Heun equations generalize many known mathematical functions including the hypergeometric ones, Mathieu functions, spheroidal wave functions, Coulomb spheroidal functions, and many others widely used in mathematical physics and applied mathematics.
It should be noted that numerous papers are devoted to expansion of solutions of the Heun equations in terms of the minor special functions; see e.g.\ \cite{El-JaickFigueiredo2013,Ishkhanyan2014,LeroyIshkhanyan2015} and references therein. 

The general Heun equation and its confluent forms appear in many fields of modern physics, such as 
general relativity, astrophysics, hydrodynamics, atomic and particle physics, etc. (see, e.g., \cite{CunhaChristiansen2011,FizievStaicova2011,RenziSammarco2016,VieiraBezerra2016,ZhangJinJinXie2016,BirkandanHortacsu2017,HartmannPortnoi2017}).
 A vast list of references to numerous physical applications, especially in general relativity, can be found in \cite{Hortacsu2017}.
``The Heun project'' (http://theheunproject.org/) should also be mentioned as a good source of information on the current development.  

Despite the increasing interest to the Heun equations, the only, to author's knowledge, software package able to evaluate the confluent Heun functions numerically is Maple{\footnotesize\texttrademark}. The purpose of the present work is to develop alternative algorithms. Following \cite{Mot15}, for numerical evaluation of the confluent Heun functions we suggest a procedure based on power series, asymptotic expansions and analytic continuation. Program realization is presented in \cite{mycode2} as Octave/Matlab code. Results of numerical tests and comparison with cases when confluent Heun functions reduce to elementary functions are given.

The proposed approach is applicable for computation of the multi-valued confluent Heun functions. We also define their single-valued counterparts by fixation of branch cuts. For the single-valued functions, an improvement of the algorithm for points close to the singular ones is suggested.

The algorithms of this work are not intended to be universal. Surely, numerical problems are expected and special treatment is needed e.g. for the cases of merging singular points (see \cite{LaySlavyanov1999}) or large accessory parameter.

\section{Statement and basic notations} 
\label{sect:statement}

We use the following form of the confluent Heun equation:
\begin{equation}
  \MF{\HeunC''}(z)+\left(\frac{\gamma}{z}+\frac{\delta}{z-1}+\varepsilon\right)\MF{\HeunC'}(z)+
  \frac{\alpha z-q}{z(z-1)}\MF{\HeunC}(z)=0.
\label{eq:HeunC}
\end{equation}
This second order linear differential equation has regular singularities at $z = 0$ and $1$, and an irregular singularity of rank 1 at $z = \infty$ (see e.g.\ \cite{SlavyanovLay2000}).
The parameter $q\in\Cb$ is usually referred to as an accessory or
auxiliary parameter and $\gamma$, $\delta$, $\varepsilon$, $\alpha$
(also belonging to $\Cb$) are exponent-related parameters.
It is important to note that in this paper the parameters $\varepsilon$ and $\alpha$ are assumed to
be independent. Below, we will use notation $\HeunC(q,\alpha,\gamma,\delta,\varepsilon;z)$ or $\HeunC(z)$ for brevity.

There are $6$ local solutions of equation \eqref{eq:HeunC}. The Frobenius method can be used to derive local power-series solutions to
\eqref{eq:HeunC} near $z=0$ and $z=1$ (two per a singular point), while two solutions at $z=\infty$ can be obtained in the form of asymptotic series.  In \S\;\ref{sect:expat0} we will present the local solutions near the point $z = 0$. One of them is
analytic in a vicinity of zero and if $\gamma$ is not a nonpositive integer, we normalize this solution to unity at $z=0$ and call it the local confluent Heun function. It is denoted by $\HeunCl(q,\alpha,\gamma,\delta,\varepsilon;z)$. For the second Frobenius local
solution, we will use the notation $\HeunCs(q,\alpha,\gamma,\delta,\varepsilon;z)$.

When $\gamma$ is a nonpositive integer, one solution of \eqref{eq:HeunC} is
analytic in a vicinity of $z=0$ but it is equal to zero at $z=0$, whereas the second solution can be
normalized to unity at zero but generally it is not analytic. Following \cite{Mot15},
the normalized solution will be denoted by $\HeunCl(z)$ and another one by $\HeunCs(z)$.

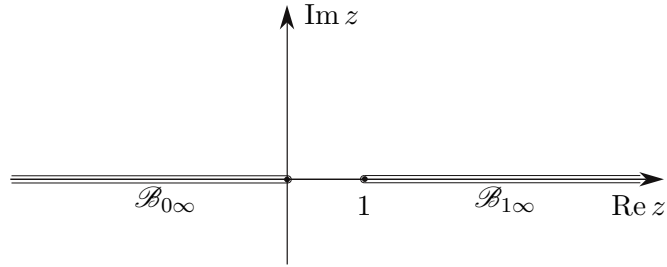
\begin{figure}[h!]
\centering
 \SetLabels
 \L (0.53*0.19) $1$\\
 \L (0.45*0.925) $\Im z$\\
 \L (0.92*0.19) $\Re z$\\
 \L (0.19*0.22) $\BCzero$\\
 \L (0.71*0.22) $\BCone$\\
 \endSetLabels
 \leavevmode\AffixLabels{\includegraphics[width=86mm]{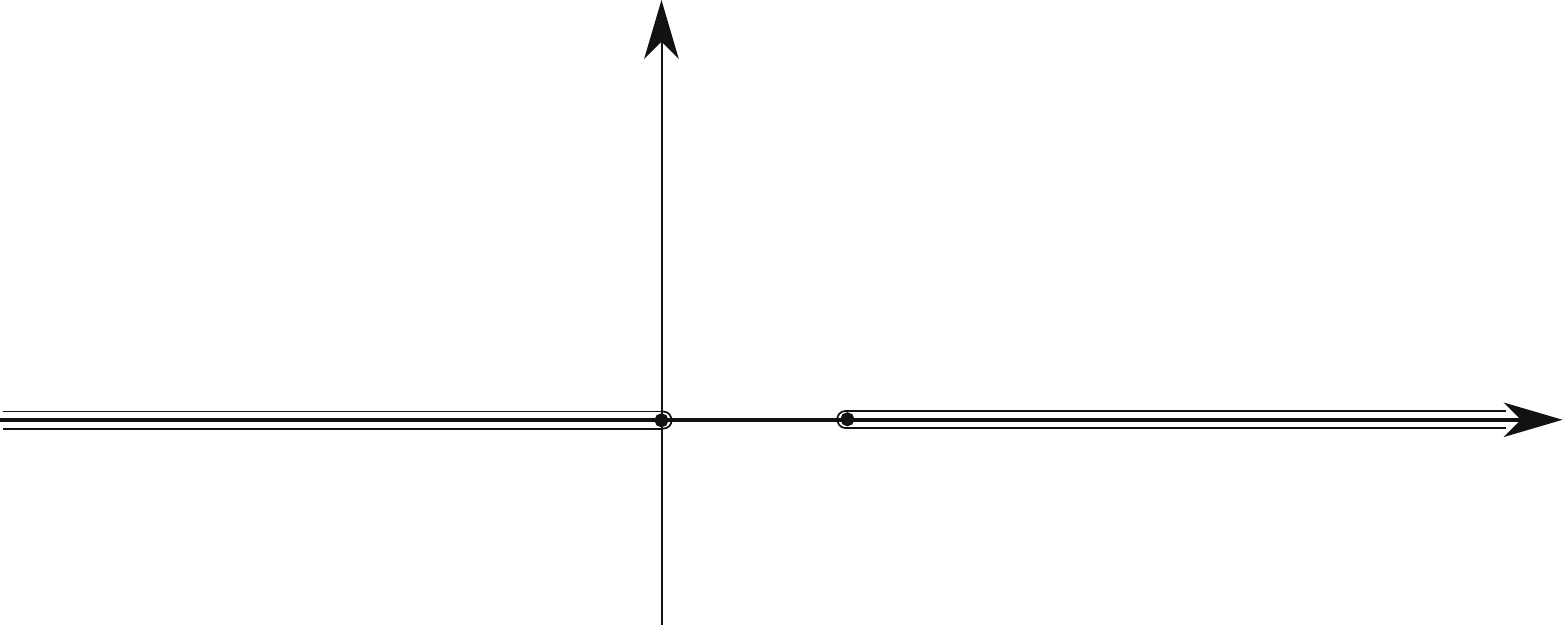}}
 \vspace{-1mm}
 \caption{Branch cuts.}
 \vspace{0mm}
 \label{bcwithoutlog}
\end{figure}

It is important to note that generally $\HeunC(z)$ is a multi-valued function and, so, to define single-valued
functions $\HeunCl(z)$ and $\HeunCs(z)$, we should choose branch cuts. In the present work, we fix the branch
cuts $\BCzero=(-\infty,0)$ and $\BCone=(1,+\infty)$, connecting the points $0$ and $1$ to~$\infty$, respectively (see Fig.~\ref{bcwithoutlog}). 
For $\gamma\not\in\Zzo$ ($\Zzo$ means the set consisting of zero and negative integers), for definition of single-valued $\HeunCl(z)$ it is sufficient to use $\BCone$. It is the case for $\HeunCs(z)$ when $\gamma\in\Zzo$.

\section{Power series expansions at the point $\bm{z=0}$}
\label{sect:expat0}

Power series expansion of the confluent Heun function $\HeunCl(z)$, such that
$\HeunCl(0)=1$, is well-known
for $\gamma\not\in\Zzo$. We have
\begin{equation}
\HeunCl(q,\alpha,\gamma,\delta,\varepsilon;z)=\sum_{n=0}^{\infty}b_n z^n,
\label{eq:HeunL0series}
\end{equation}
where the coefficients $b_n$ are submitted to the following three-term
recurrence relation:
\begin{equation}
  P_n b_n=Q_n b_{n-1}+R_n b_{n-2}.
\label{eq:HeunL0Andef}
\end{equation}
Here
\begin{equation}
\begin{gathered}
P_n = n (\gamma-1+n),\qquad
Q_n = -q+(n-1)(\gamma+\delta-\varepsilon+n-2),\qquad 
R_n = (n-2)\varepsilon+\alpha,
\end{gathered}
\label{eq:PQRdef}
\end{equation}
and the initial conditions are as follows: $b_{-1}=0$, $b_0=1$. (Then $\HeunCl'_z(q,\alpha,\gamma,\delta,\varepsilon;0)=-q/\gamma$.)

The confluent Heun function
$\HeunCl(z)$ is analytic in the circle $|z|<1$ and
Cauchy's theorem on the expansion of an analytic function into a power series
(see e.g.\ Theorem 16.7 in \cite[Part~I]{Markushevich1977}) guarantees that the
series \eqref{eq:HeunL0series} converges to $\HeunCl(z)$ inside the circle
$|z|<1$. (Though the question of the forward stability of the recursion relations, see e.g.\ \cite{Wimp1984}, is out of our scope in this paper.) 

In the case of integer $\gamma$, the local Frobenius solution corresponding to the
smaller exponent ($0$ or $1-\gamma$) may contain a logarithmic factor (see e.g.\
\cite{Kamke1944,Wasow}). So, for $\gamma\in\Zzo$ we are looking for the
solution of \eqref{eq:HeunC} in the following form:
\begin{equation}
\HeunCl(q,\alpha,\gamma,\delta,\varepsilon;z)=\sum_{n=0,\,n\neq\nsp}^{\infty}\!\!c_n
z^n +\log(z)\sum_{n=\nsp}^{\infty}s_n z^n,
\label{eq:HeunLserieslog}
\end{equation}
where $\nsp=1-\gamma$. Note that a solution with the sought property
$\HeunCl(0)=1$ could be found for any $c_{\nsp}$. We fix $c_{\nsp}=0$ for
definiteness.

Substituting \eqref{eq:HeunLserieslog} into \eqref{eq:HeunC}, we find
\begin{equation}
\log(z)\MF{\Lkr}\Biggl(\,\sum_{n=\nsp}^\infty s_n z^n\Biggr)+
\MF{\Lkr}\Biggl(\,\sum_{n=0,\,n\neq\nsp}^\infty\! c_n z^n\Biggr)+
\MF{\Hat{\Lkr}}\Biggl(\,\sum_{n=\nsp}^\infty s_n z^n\Biggr)=0,
\label{eq:HeunEqLserieslog}
\end{equation}
where $\Lkr$ is the operator of the confluent Heun equation such that \eqref{eq:HeunC} is
written as $\Lkr\!\HeunC=0$. Besides,
\[
\left(\Hat{\Lkr}\psi\right)(z)=\frac2z\MF{\psi'}(z)+\frac1z\left(
\frac{\gamma-1}{z}+\frac{\delta}{z-1}+\varepsilon\right)\MF{\psi}(z).
\]
Obviously, \eqref{eq:HeunEqLserieslog} splits into two equations
\begin{gather}
\MF{\Lkr}\Biggl(\,\sum_{n=\nsp}^\infty s_n
z^n\Biggr)=0,\label{eq:HeunEqLserieslog1}\\
\MF{\Lkr}\Biggl(\,\sum_{n=0,\,n\neq\nsp}^\infty\! c_n z^n\Biggr)+
\MF{\Hat{\Lkr}}\Biggl(\,\sum_{n=\nsp}^\infty s_n z^n\Biggr)=0,
\label{eq:HeunEqLserieslog2}
\end{gather}

Let us collect in \eqref{eq:HeunEqLserieslog1}, \eqref{eq:HeunEqLserieslog2} terms having the same asymptotic
nature as $z\to0$. First, we find that coefficients $c_n$ for $n=1,\ldots,\nsp-1$
are submitted to the recurrence \eqref{eq:HeunL0Andef}: $P_n c_n=Q_n c_{n-1}+R_n
c_{n-2}$, where $P_n$, $Q_n$, $R_n$ are defined by \eqref{eq:PQRdef} and the
initial conditions are $c_{-1}=0$, $c_0=1$.

From \eqref{eq:HeunEqLserieslog1} and
\eqref{eq:HeunEqLserieslog2}, we find that the coefficients $s_n$ for
$n=\nsp+1,\nsp+2,\ldots$ are submitted to the same recurrence relationship
\eqref{eq:HeunL0Andef}: $P_n s_n=Q_n s_{n-1}+R_n s_{n-2}$, where $s_{\nsp-1}=0$
and another initial condition includes coefficients $c_{\nsp-1}$, $c_{\nsp-2}$:
\begin{gather*}
\nsp\, s_{\nsp}=c_{\nsp-1}\bigl[-q+\gamma(1-\delta+\varepsilon)\bigr]+
c_{\nsp-2}\bigl[\alpha-\varepsilon(1+\gamma)\bigr].
\end{gather*}

At the next step, we can define coefficients $c_n$ for $n=\nsp+1,\nsp+2,\ldots$\,\
From \eqref{eq:HeunEqLserieslog2} we obtain the following relationship:
\begin{equation}
P_n c_n=Q_n c_{n-1}+R_n c_{n-2}+ S_n s_{n}+T_n s_{n-1}+U_n s_{n-2},
\label{eq:recurrlog}
\end{equation}
where
\begin{gather*}
S_n = 1-\gamma-2n,\qquad
T_n = \gamma+\delta-\varepsilon+2n-3,\qquad
U_n = \varepsilon.
\end{gather*}

In this way for $\gamma\in\Zzo$, using \eqref{eq:HeunLserieslog},  we
obtain a local solution, equal to unity at $z=0$. 
For the constructed solution, it is easy to find that 
$\HeunCl'(q,\alpha,\gamma,\delta,\varepsilon;0)=-q/\gamma$  for $\gamma\in\Zzo\setminus\{0\}$. For $\gamma=0$, we have
$\HeunCl'(q,\alpha,\gamma,\delta,\varepsilon;0)/\log(z)\to-q$ as $z\to0$.

\begin{samepage}
The second
local solution can be defined as follows (see also
\eqref{eq:HeunEqLserieslog1}):
\begin{equation}
  \HeunCs(q,\alpha,\gamma,\delta,\varepsilon;z)=\sum_{n=\nsp}^{\infty}\check{s}_n z^n,
\label{eq:HeunSsergammanonpositive}
\end{equation}
where $P_n \check{s}_n=Q_n \check{s}_{n-1}+R_n \check{s}_{n-2}$ for $n>\nsp$ and
$\check{s}_{\nsp}=1$, $\check{s}_{\nsp-1}=0$. 

\end{samepage}

As it was mentioned above, $c_{\nsp}$ in \eqref{eq:HeunLserieslog} could be arbitrary. In other
words, the above choice of $\HeunCl(q,\alpha,\gamma,\delta,\varepsilon;z)$ for
$\gamma\in\Zzo$ is non-unique; it could be a linear combination
\begin{equation*}
  \HeunCl(q,\alpha,\gamma,\delta,\varepsilon;z)+C\HeunCs(q,\alpha,\gamma,\delta,\varepsilon;z)
\label{eq:lincomb}
\end{equation*}
for an arbitrary constant $C$.

For $\gamma\not\in\Zzo$, by substitution to \eqref{eq:HeunC} it is straightforward to check the following relationship:
\begin{equation}
\HeunC(q,\alpha,\gamma,\delta,\varepsilon;z)=\E{-\varepsilon z}\HeunC(q-\varepsilon\gamma,\alpha-\varepsilon(\gamma+\delta),\gamma,\delta,-\varepsilon;z).
\label{eq:expout}
\end{equation}
(The latter formula is pretty useful in numerical evaluation for large $|\varepsilon z|$.) However, for $\gamma\in\Zzo$ (due to the non-uniqueness) the formula \eqref{eq:expout} is generally not true for $\HeunCl(z)$. Namely, we have 
\begin{equation}
\HeunCl(q,\alpha,\gamma,\delta,\varepsilon;z)+\mathcal{A}\HeunCs(q,\alpha,\gamma,\delta,\varepsilon;z)=\E{-\varepsilon z}\HeunCl(q-\varepsilon\gamma,\alpha-\varepsilon(\gamma+\delta),\gamma,\delta,-\varepsilon;z),
\label{eq:expout-mod}
\end{equation}
where
\[
\mathcal{A}=-\sum_{n=0}^{\nsp} c_n \frac{\varepsilon^{\nsp-n}}{(\nsp-n)!}.
\]

Consider now the function $\HeunCs(z)$ for arbitrary $\gamma$. We
should discern two situations: $\gamma=1$ and $\gamma\neq1$. In the latter case, we can use the following representation:
\begin{equation}
 \HeunCs(q,\alpha,\gamma,\delta,\varepsilon;z) = z^{1-\gamma}
 \HeunCl\bigl(q+(\gamma-1)(\delta-\varepsilon),\alpha+\varepsilon(1-\gamma),2-\gamma,\delta,\varepsilon;z\bigr).
\label{eq:HeunSgammaneq1}
\end{equation}
Notably, this formula includes \eqref{eq:HeunSsergammanonpositive} as a
particular case, justifying our way to introduce $\HeunCl(z)$ and $\HeunCs(z)$ for
non-positive integer $\gamma$. For $\gamma=2,3,\ldots$,  the expansion of the function $\HeunCs(z)$ defined by \eqref{eq:HeunSgammaneq1} and \eqref{eq:HeunLserieslog} contains logarithmic term.

For $\gamma=1$, repeating the arguments used to derive representation of
$\HeunCl(z)$ in the case $\gamma\in\Zzo$, we can find the following local
representation
\begin{equation}
\HeunCs(q,\alpha,\gamma,\delta,\varepsilon;z)=\sum_{n=1}^{\infty}d_n z^n
+\log(z)\sum_{n=0}^{\infty}t_n z^n.
\label{eq:HeunSgamma=1ser}
\end{equation}
Here $P_n t_n=Q_n t_{n-1}+R_n t_{n-2}$, where $t_{-1}=0$, $t_{0}=1$ and (cf.\
\eqref{eq:recurrlog})
\begin{equation*}
P_n d_n=Q_n d_{n-1}+R_n d_{n-2}+ S_n t_{n}+T_n t_{n-1}+U_n t_{n-2},\quad
d_{-1}=d_0=0.
\end{equation*}

\section{Expansions at infinity}
\label{sect:expatinf}

In this section we write expansions of the confluent Heun function at infinity, where the equation has an irregular singularity of rank 1.
Assuming that $\varepsilon\neq0$, we look for a solution in the form
\begin{equation}
\HeunC_{\!A,\infty}(q,\alpha,\gamma,\delta,\varepsilon;z)=(-z)^{-\frac{\alpha}{\varepsilon}}\sum_{n=0}^{\infty} \beta_n \frac{n!}{(\varepsilon z)^n}.
\label{eq:HeunAinfseries}
\end{equation}
From \eqref{eq:HeunC}, we find that the coefficients $\beta_n$ are subject to the 
recurrence
\begin{equation*}
  \beta_n=\tilde{Q}_n \beta_{n-1}+\tilde{R}_n \beta_{n-2},
\label{eq:HeunAinfdef}
\end{equation*}
where
\begin{equation*}
\begin{gathered}
\tilde{Q}_n = 1+n^{-2}\biggl(-q+\frac{\alpha}{\varepsilon}\biggl(2n-\gamma-\delta-1+\frac{\alpha}{\varepsilon}\biggr)
 +(\gamma+\delta-\varepsilon+1)(1-n)+\alpha-1\biggr),\\ 
\tilde{R}_n = \varepsilon\frac{\Bigl(n-2+\frac{\alpha}{\varepsilon}\Bigr)\Bigl(\gamma-n+1-\frac{\alpha}{\varepsilon}\Bigr)}{n^2(n-1)},
\end{gathered}
\label{eq:AinfPQRdef}
\end{equation*}
and the initial conditions are chosen to be as follows:
\[
\beta_{-1}=0,\qquad \beta_0=1.
\]

The second solution can be introduced by using the relationship \eqref{eq:expout}:
\begin{equation}
\HeunC_{\!B,\infty}(q,\alpha,\gamma,\delta,\varepsilon;z)=\E{-\varepsilon z}\HeunC_{\!A,\infty}(q-\varepsilon\gamma,\alpha-\varepsilon(\gamma+\delta),\gamma,\delta,-\varepsilon;z).
\label{eq:expoutinf}
\end{equation}

We note that in view of \eqref{eq:expoutinf} the so-called Stokes (anti-Stokes) line can be defined as $\Im\{\varepsilon z\}=0$ ($\Re\{\varepsilon z\}=0$).

In the special case $\varepsilon=0$, $\alpha\neq0$, one can find two solutions in the form of the following asymptotic series
\begin{equation*}
\HeunC_\pm(z)=z^\Lambda\exp\bigl\{\pm 2\ii \sqrt{\alpha z}\bigr\}\sum_{n=0}^{\infty}\beta^\pm_n z^{-\frac{n}{2}},
\label{eq:infeps0}
\end{equation*}
where
\[
\Lambda=\frac14-\frac{\gamma+\delta}{2}
\]
and the coefficients $\beta^\pm_n$ are defined by the recursion:
\begin{equation*}
\mathring{P}_n\beta^\pm_n=\pm\mathring{Q}_n \beta^\pm_{n-1}+
\mathring{R}_n \beta^\pm_{n-2}\pm\mathring{S}_n \beta^\pm_{n-3},
\label{eq:infeps0rec}
\end{equation*}
with the initial conditions $\beta^\pm_{-2}=\beta^\pm_{-1}=0$, $\beta^\pm_{0}=1$. 
In the recurrence relation,
\begin{gather*}
\mathring{P}_n=4\ii n\sqrt{\alpha},\qquad
\mathring{Q}_n=\left(n-\tfrac32\right)\left(n+\tfrac12\right)+4(\alpha-q)-(\gamma+\delta)(\gamma+\delta-2),\\
\mathring{R}_n=4\ii\sqrt{\alpha}\left(n-2+\delta\right),\quad
\mathring{S}_n= -\left(n-\tfrac32-\gamma+\delta\right)\left(n-\tfrac72+\gamma+\delta\right).
\end{gather*}
 
 
In the case $\alpha=\varepsilon=0$, the confluent Heun equation reduces to the hypergeometric one (see \cite{Erdelyi1955}, Ch.\;2) and the point $z=\infty$ is regular singular. 

\section{Power series expansion at an arbitrary regular point}
\label{sect:expatarbpt}

Further we will extend the local confluent Heun functions outside the circle of
convergence of the series \eqref{eq:HeunL0series}, \eqref{eq:HeunLserieslog},
\eqref{eq:HeunSgamma=1ser} ($|z|<1$). For this purpose, in \S~\ref{sect:basic_algorithm} we will use analytic continuation process based on the power series expansion which we derive in this section.

We seek the solution $\HeunC_{\!(z_0,H_0,H'_0)}(q,\alpha,\gamma,\delta,\varepsilon;z)$
to equation \eqref{eq:HeunC} satisfying the conditions
\begin{equation}
 \HeunC(q,\alpha,\gamma,\delta,\varepsilon;z_0)=H_0,\qquad
 \frac{\partial}{\partial z}\HeunC(q,\alpha,\gamma,\delta,\varepsilon;z)\Bigr|_{z=z_0}=H'_0.
\label{eq:z0H0H'0}
\end{equation}
Here $z_0$ is an arbitrary finite point, assumed not to coincide with the singular points
$0$, $1$. We look for power series expansion of the confluent Heun function in
the form
\begin{equation}
\HeunC_{\!(z_0,H_0,H'_0)}(q,\alpha,\gamma,\delta,\varepsilon;z)=\sum_{n=0}^{\infty}\car_n(z-z_0)^n.
\label{eq:HeunZ0series}
\end{equation}

Substituting \eqref{eq:HeunZ0series} into \eqref{eq:HeunC} and collecting terms
at the same power of $z-z_0$, we obtain the following 4-term recurrence
relation for the coefficients $\car_n$:
\begin{equation}
\Pcoef_n \car_n=\Qcoef_n \car_{n-1}+\Rcoef_n \car_{n-2} + \Scoef_n
\car_{n-3},\label{eq:HeunZ0Cndef}
\end{equation}
where
\begin{gather*}
\Pcoef_n=n(1-n)z_0(z_0-1),\qquad
\Qcoef_n=(n-1)
\bigl(\varepsilon z_0^2+z_0(\gamma+\delta-\varepsilon+2(n-2))-\gamma-n+2\bigr), \\
\Rcoef_n=z_0(2(n-2)\varepsilon+\alpha)+(n-2)(\gamma-\varepsilon+\delta+n-3)-q,\qquad
\Scoef_n=(n-3)\varepsilon+\alpha.
\end{gather*}
Obviously, the solution \eqref{eq:HeunZ0series} satisfies the conditions \eqref{eq:z0H0H'0} if
the recurrence process starts with the initial conditions 
\begin{equation*}
\car_{-1}=0,\qquad \car_{0}=H_0,\qquad \car_{1}=H'_0.
\end{equation*}

The series \eqref{eq:HeunZ0series} converges inside the circle $|z-z_0|<r$,
where $r$ is the distance to the nearest singular point,
$r=\min\{|z_0|,|z_0-1|\}$. Of course, practically the convergence can be slow when $|z-z_0|$ is not small.

\section{Basic algorithm}
\label{sect:basic_algorithm}

Let us introduce the projection operator $\mathscr{P}_{z}^N$ which, being applied to
an analytic function, truncates its power series expansion at the point $z$ to
the first $N$ terms. Consider first $\gamma\not\in\Zzo$. Using the expansion \eqref{eq:HeunL0series}, we evaluate
\begin{equation}
\bigl(\mathscr{P}_0^N\!\HeunCl\bigr)(z)=\sum_{n=0}^{N}b_n z^n,\qquad
\bigl(\mathscr{P}_0^N\!\HeunCl\bigr)'(z)=\sum_{n=1}^{N}n\, b_n z^{n-1},
\label{eq:HeunLpwrserappr}
\end{equation}
as  approximation of $\HeunCl(z)$ and $\HeunCl'(z)$  in a vicinity of $z=0$.

In our algorithm we do not fix the number $N$ in the representations
\eqref{eq:HeunLpwrserappr}; it will be defined as we proceed with
computation of series terms and summation until a termination
condition is satisfied. Namely, we stop the process when
$\bigl(\mathscr{P}_0^{N}\!\HeunCl\bigr)(z)$,
$\bigl(\mathscr{P}_0^{N-1}\!\HeunCl\bigr)(z)$, and
$\bigl(\mathscr{P}_0^{N}\!\HeunCl\bigr)'(z)$,
$\bigl(\mathscr{P}_0^{N-1}\!\HeunCl\bigr)'(z)$ are not distinguishable in the
used computer arithmetics.

To estimate the quality of the approximation, in view of \eqref{eq:HeunC} we
compute the value
\[
  \underline{\HeunCl}(z)=\frac{1}{q-\alpha z}\biggl\{z(z-1)\bigl(\mathscr{P}_0^N\!\HeunCl\bigr)''(z)+
  \bigl[\gamma(z-1)+\delta z+\varepsilon z(z-1)\bigr]\bigl(\mathscr{P}_0^N\!\HeunCl\bigr)'(z)\biggr\},
\]
where $\bigl(\mathscr{P}_0^N\!\HeunCl\bigr)''(z)=\sum_{n=2}^{N}n(n-1)\,b_n z^{n-2}$.
Then we suppose proximity of
\begin{equation}
\MF{r_0}(z)=\bigl|\underline{\HeunCl}(z)-\bigl(\mathscr{P}_0^N\!\HeunCl\bigr)(z)\bigr|
\label{eq:errdef}
\end{equation}
to the true error of the approximation
$\bigl|\HeunCl(z)-\bigl(\mathscr{P}_0^N\!\HeunCl\bigr)(z)\bigr|$.  Near the point
$z=z_*=q/\alpha$, numerical computation of $\underline{\HeunCl}(z)$ is unreliable due to essential loss of significance.  In a vicinity of $z_*$, it
can be suggested  to use an estimate based on properties of the series, e.g., akin to one used in \cite{mycode2}, 
\begin{equation}
\MF{\Hat{r}\vphantom{r}_0}(z)=
\sqrt{N}\bigl|z^Nb_N(z)\bigr| + \epsilon N \bigl|\bigl(\mathscr{P}_0^N\!\HeunCl\bigr)(z)\bigr|,
\label{eq:err2def}
\end{equation}
where $\epsilon$ is machine epsilon in the applied computer arithmetics.

We write the described algorithm as a function $\psHeunCl(z)$ which returns 4-tuple
\[
 \psHeunCl: z\mapsto [f,f',r,N],
\]
where $N+1$ is the number of terms in power series, defined by the termination
condition, $f=\bigl(\mathscr{P}_0^N\!\HeunCl\bigr)(z)$,
$f'=\bigl(\mathscr{P}_0^N\!\HeunCl\bigr)'(z)$, and $r$ is the value computed with \eqref{eq:errdef} or
\eqref{eq:err2def}.

The scheme of computation of $\HeunCl(z)$ in the case $\gamma\in\Zzo$
is analogous, but slightly more involved. We use
\eqref{eq:HeunLserieslog} and, instead of \eqref{eq:HeunLpwrserappr}, define the
function $\psHeunCl(z)$ starting from the expression
\[
\sum_{n=0,\,n\neq\nsp}^{N}\!\!c_n z^n +\log(z)\sum_{n=\nsp}^{N}s_n z^n.
\]

Assume that $|z|<\reservecoeff$, where  $\reservecoeff\in(0,1)$ is some coefficient chosen
so that $N$ defined by the termination condition is expected to be moderate (in computations presented in \S~\ref{sect:num}, we fix $\reservecoeff=0.38$). Then we can use the numerical algorithm $\psHeunCl$ for
evaluation of the function $\HeunCl(z)$ and its derivative, and for estimation of
the approximation error.

Consider further the case $|z|\geq\reservecoeff$. First we define an
auxiliary algorithm. Let $z_0$ be an arbitrary point not belonging to the set
$\{0,1,\infty\}$. Using \eqref{eq:HeunZ0series} we define
\[
\bigl(\mathscr{P}_{\!z_0}^N\HeunC_{\!(z_0,H_0,H'_0)}\bigr)(z)=\sum_{n=0}^{N}\car_n(z-z_0)^n,\qquad
\bigl(\mathscr{P}_{\!z_0}^N\HeunC_{\!(z_0,H_0,H'_0)}\bigr)'(z)=\sum_{n=1}^{N}n\,\car_n(z-z_0)^{n-1},
\]
as approximations of $\HeunC_{\!(z_0,H_0,H'_0)}(z)$ and
$\HeunC'_{(z_0,H_0,H'_0)}(z)$ for $z$ close to $z_0$. Here coefficients $\car_n$ are defined by  \eqref{eq:HeunZ0Cndef} and we proceed with summation until the termination condition (analogous to that described above) is satisfied.

Again, we compute
\begin{multline*}
  \underline{\HeunC\!}\vphantom{\HeunC}_{\,(z_0,H_0,H'_0)}(z)=\frac{1}{q-\alpha z}\
  \biggl\{z(z-1)\bigl(\mathscr{P}_{\!z_0}^N\HeunC_{\!(z_0,H_0,H'_0)}\bigr)''(z)\\
  {}+
  \Bigl[\gamma(z-1)+\delta z+\varepsilon z(z-1)\Bigr]\bigl(\mathscr{P}_{\!z_0}^N
  \HeunC_{\!(z_0,H_0,H'_0)}\bigr)'(z)\biggr\},
\end{multline*}
and the value
$\MF{r_{(z_0,H_0,H'_0)}}(z)=\bigl|\underline{\HeunC\!}\vphantom{\HeunC}_{\,(z_0,H_0,H'_0)}(z)-
\bigl(\mathscr{P}_{\!z_0}^N\HeunC_{\!(z_0,H_0,H'_0)}\bigr)(z)\bigr|$. In view of
essential loss of significance in computation of
$\underline{\HeunC\!}\vphantom{\HeunC}_{\,(z_0,H_0,H'_0)}(z)$ near $z=z_*$, we define
\begin{equation*}
\MF{\Hat{r}\vphantom{r}_{(z_0,H_0,H'_0)}}(z)=
\sqrt{N}\bigl|(z-z_0)^N\car_N(z)\bigr| + \epsilon N \bigl|\bigl(\mathscr{P}_{\!z_0}^N\HeunC_{\!(z_0,H_0,H'_0)}\bigr)(z)\bigr|.
\end{equation*}

We write the described algorithm as a function 
$\psHeunCo_{\!(z_0,H_0,H'_0)}: z\mapsto [f,f',r,N]$,
where $N+1$ is the number of terms in power series defined by the termination
condition, $f=\bigl(\mathscr{P}_{z_0}^N\HeunC_{\!(z_0,H_0,H'_0)}\bigr)(z)$,
$f'=\bigl(\mathscr{P}_{z_0}^N\HeunC_{\!(z_0,H_0,H'_0)}\bigr)'(z)$, $r$ is equal to $\MF{r_{(z_0,H_0,H'_0)}}(z)$ or $\MF{\Hat{r}\vphantom{r}_{(z_0,H_0,H'_0)}}(z)$.

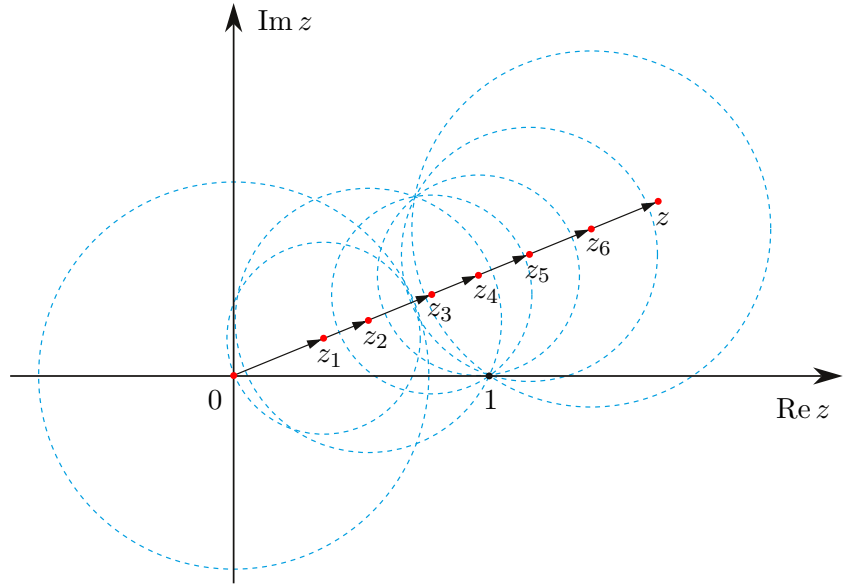
\begin{figure}[t]
\centering
 \SetLabels
 \L (0.57*0.305) $1$\\
 \L (0.3*0.95) $\Im z$\\
 \L (0.92*0.29) $\Re z$\\
 \L (0.24*0.305) $0$\\
 \L (0.775*0.619) $z$\\
 \L (0.37*0.385) $z_1$\\
 \L (0.425*0.416) $z_2$\\
 \L (0.502*0.463) $z_3$\\
 \L (0.557*0.494) $z_4$\\
 \L (0.619*0.529) $z_5$\\
 \L (0.692*0.575) $z_6$\\
 \endSetLabels
 \leavevmode\AffixLabels{\includegraphics[width=110mm]{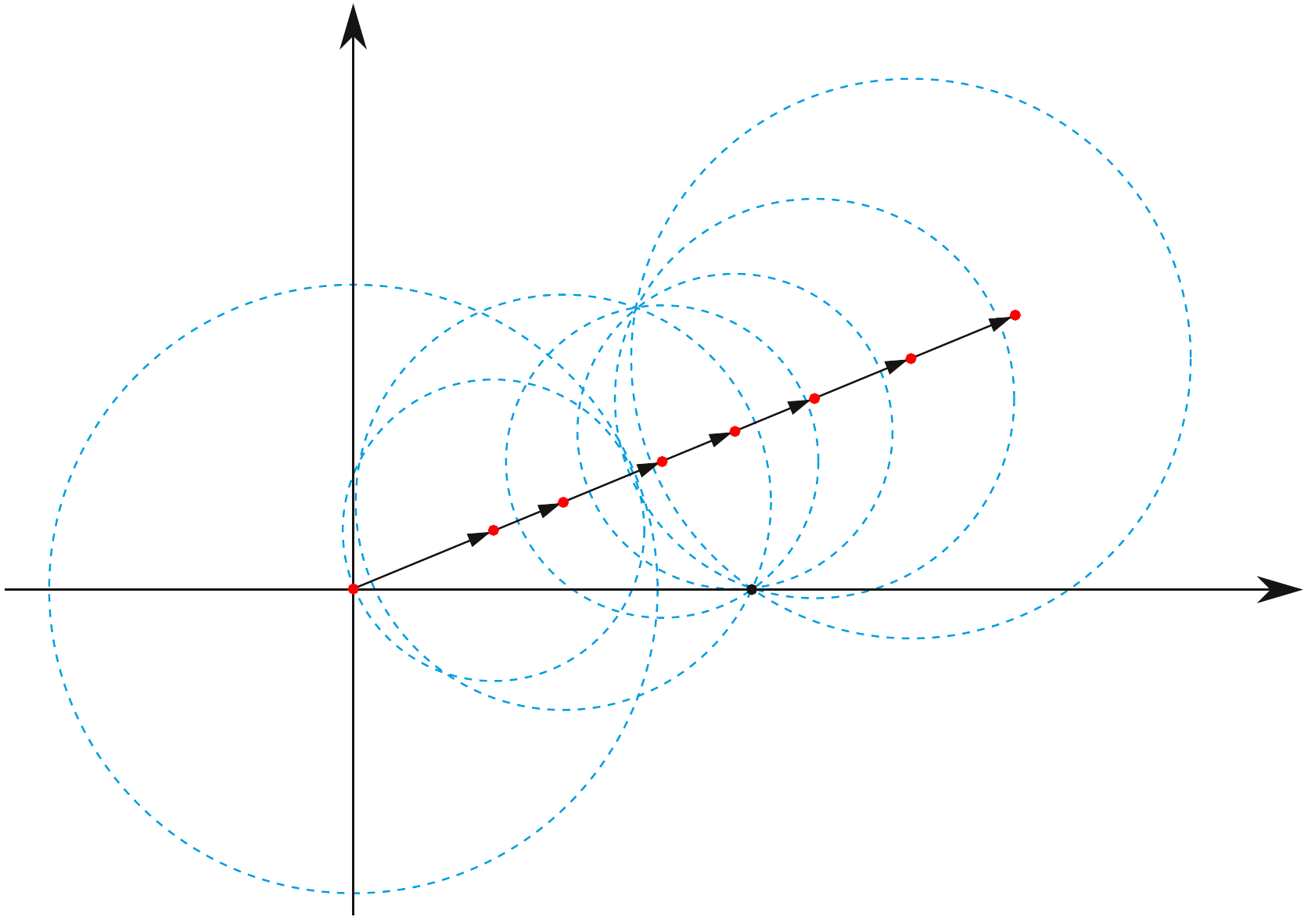}}
 \vspace{-1mm}
 \caption{Analytic continuation using power series.}
 \vspace{0mm}
 \label{analytcont}
\end{figure}


Now we are ready to proceed with analytic continuation along a path from zero to $z$. Consider first the simplest case when the path is the line segment $(0,z)$. At the first step,
we compute
\[
 [\HeunCl_1,\HeunCl'_1,r_1,N_1]=\psHeunCl(z_1),
\]
where $z_1=\E{\ii\arg(z)}\reservecoeff$ (see Fig.~\ref{analytcont}, where for definiteness we choose $\reservecoeff=0.5$). 
Further, we connect two regular points $z_1$ and $z$, starting with the values $\HeunCl_1$, $\HeunCl'_1$ at $z_1$. Denote this algorithm by $\psHeunC_{\!(z_1,\HeunCltxst_1, \HeunCltxst'_1)}$. For $p=1$, $2$, and so on, we define 
\[
R_p=\reservecoeff\min\{|z_p|,|z_p-1|\},
\]
\[
z_{p+1}=\begin{cases}z_p+\E{\ii\arg(z)} R_p & \mbox{if}\ \ |z-z_p|> R_p,\\
z & \mbox{if}\ \ |z-z_p|\leq R_p,
\end{cases}
\]
and compute
\[
 [\HeunCl_{p+1},\HeunCl'_{p+1},r_{p+1},N_{p+1}]=\psHeunCo_{(z_p,\HeunCltxst_p,\HeunCltxst'_p)}(z_{p+1}).
\]
The iterations stops when $z_{p+1}=z$. Finally, we have $\HeunCl_{p+1}$,
$\HeunCl'_{p+1}$ as approximations of $\HeunCl(z)$ and $\HeunCl'(z)$, respectively.
We also compute the values $r_\Sigma=r_1+\ldots+r_{p+1}$ and
$N_\Sigma=N_1+\ldots+N_{p+1}+p+1$. Here $N_\Sigma$ is the total number of power
series terms which can be used as a measure of computer load and $r_\Sigma$ may
be an indicator of the approximation quality.

It can be useful to modify this algorithm by allowing more precise selection of $R_p$. For example, in \cite{mycode2}, after a step of iteration is complete, we choose $R=R_p N_\diamond/N_p$, and at the next step $R_{p+1}=\min\bigl\{R,\reservecoeff\min\{|z_{p+1}|,|z_{p+1}-1|\}\bigr\}$. Here $N_\diamond$ is a number of series terms which is considered as in some sense optimal for the used computer arithmetics (in the computations of \S~\ref{sect:num}, $N_\diamond=40$).


We also note that in view of \eqref{eq:expoutinf}, if $\Re(-\varepsilon z)>0$, it may reasonable to compute $\HeunCl(q,\alpha,\gamma,\delta,\varepsilon;z)$ through $\HeunCl(q-\varepsilon\gamma,\alpha-\varepsilon(\gamma+\delta),\gamma,\delta,-\varepsilon;z)$; see \eqref{eq:expout}, \eqref{eq:expout-mod}. This trick is used in the code \cite{mycode2}.

The described algorithm of continuation along a line segment is readily
generalized for the case when $0$ and $z$ are connected by a polyline
$\Upsilon$. This gives us a way to compute the multi-valued confluent Heun function. The
resulting procedure can be considered as a function
$\psHeunClpath{\Upsilon}: z\mapsto [f,f',r_\Sigma,N_\Sigma]$,
where $f$ and $f'$ are the resulting approximations of the confluent Heun function at $z$
and its derivative. 

The above arguments can be literally exploited to define the function $\psHeunCspath{\Upsilon}$. The procedure of
analytic continuation described above for $\HeunCl(z)$ can be applied for $\HeunCs(z)$ with simple
modification\,---\,it should start from another expansion at $z=0$, given by \eqref{eq:HeunL0series}, \eqref{eq:HeunSgammaneq1} or by \eqref{eq:HeunSgamma=1ser}.

It is notable that the size of the step in the described analytic
continuation is small for parts of the polyline $\Upsilon$ close to a
singular point. This also means an increase of the number of used circular
elements in the continuation procedure which, in its turn, may lead to loss of
accuracy. The influence of the singular points can be reduced by a choice of the
path of continuation. 

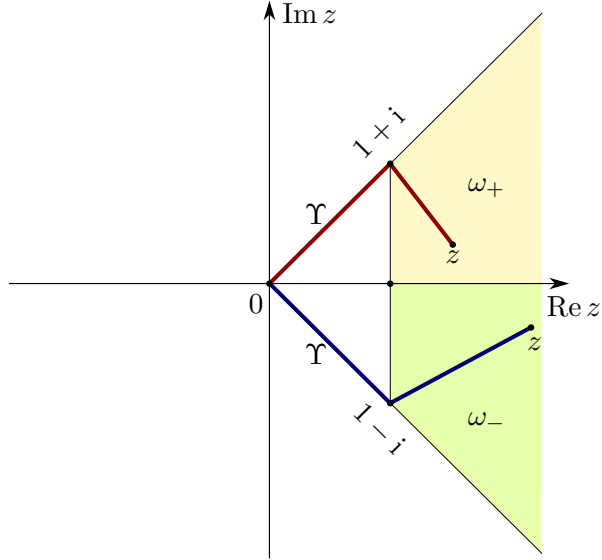
\begin{figure}[t]
\centering
 \SetLabels
 \L (0.49*0.96) $\Im z$\\
 \L (0.96*0.435) $\Re z$\\
 \L (0.61*0.715) \rotatebox{45}{$1+\ii$}\\
 \L (0.61*0.26) \rotatebox{-45}{$1-\ii$}\\
 \L (0.43*0.44) $0$\\
 \L (0.78*0.53) $z$\\
 \L (0.925*0.375) $z$\\
 \L (0.53*0.6) $\Upsilon$\\
 \L (0.53*0.35) $\Upsilon$\\
 \L (0.82*0.66) $\omega_+$\\
 \L (0.82*0.24) $\omega_-$\\
 \endSetLabels
 \leavevmode\AffixLabels{\includegraphics[width=74mm]{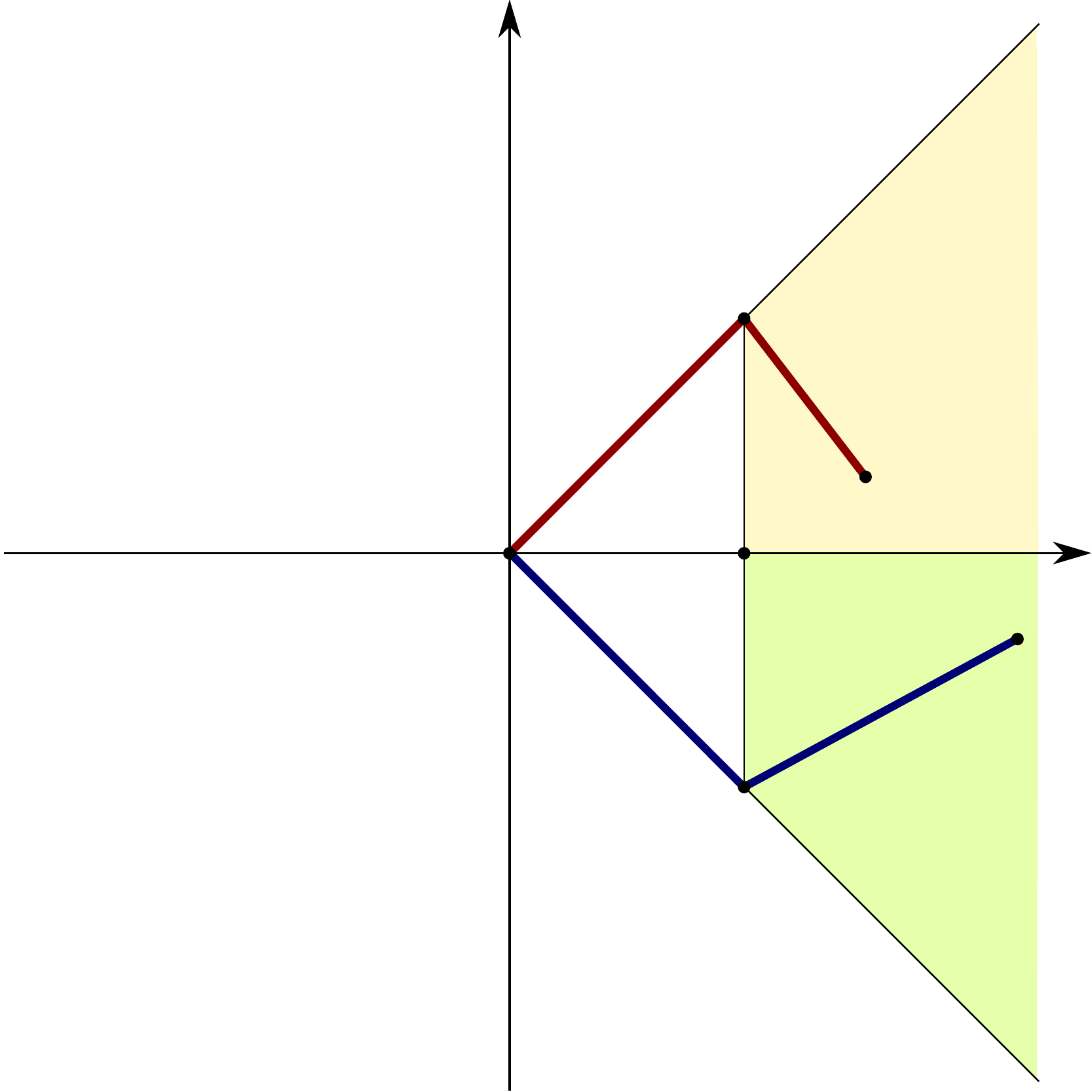}}
 \vspace{-1mm}
 \caption{Path from zero to $z$ consisting of two line segments for $z\in\omega_+$ and $z\in\omega_-$.}
 \vspace{0mm}
 \label{analytcont2}
\end{figure}

In the computational scheme applied in \S~\ref{sect:num} \cite{mycode2}, we use paths consisting of two line segments when $z$ belongs to one of the domains
$\omega_\pm=\{z:\Re z>1,0<\pm\Im z<\Re z\}$ (see Fig.~\ref{analytcont2}). 
Thus, for $z\in\omega_\pm$ we consider the path $\Upsilon$ that consequently connects the points $0$, $1+\ii\sign(\Im z)$, $z$,  and define
$\psHeunClfin(z)=\psHeunClpath{\Upsilon}(z)$ and
$\psHeunCsfin(z)=\psHeunCspath{\Upsilon}(z)$.

\section{Computation of single-valued confluent Heun functions near singular points}
\label{sect:improvement}

As it is already noted, the number of circular elements in the continuation
procedure for computation of $\HeunCl(z)$, $\HeunCs(z)$  increases as $z$ approaches a singular point
($1$ or $\infty$). In this section we suggest improvements of the algorithm near these points.

\begin{samepage}
Consider first a vicinity of $z=1$. It is known that two local solutions can be written as follows:
$\HeunCl(q-\alpha,-\alpha,\delta,\gamma,-\varepsilon;1-z)$,
$\HeunCs(q-\alpha,-\alpha,\delta,\gamma,-\varepsilon;1-z)$.
Thus, we have
\begin{equation}
\HeunCl(q,\alpha,\gamma,\delta,\varepsilon;z) = 
C_1\HeunCl(q-\alpha,-\alpha,\delta,\gamma,-\varepsilon;1-z)+
C_2\HeunCs(q-\alpha,-\alpha,\delta,\gamma,-\varepsilon;1-z),
\label{eq:near1}
\end{equation}
where $C_1$, $C_2$ are some constants.

\end{samepage}

To the author knowledge, an explicit solution to the
two-point connection problem for the confluent Heun equation has not been found (see e.g.\ \cite{Kazakov2006} and references therein). So, we define the matching
coefficients $C_1$, $C_2$ numerically, in the following way. We choose a matching point,
$z=\midpoint{1}=1/2$ and apply the algorithms $\psHeunClfin$ and $\psHeunCsfin$ described in
\S\;\ref{sect:basic_algorithm} to find
\begin{gather*}
 [f_0,f'_0,r_0,N_0]=\psHeunClfin(q,\alpha,\gamma,\delta,\varepsilon;\midpoint{1}),\\
 [f_1,f'_1,r_1,N_1]=\psHeunClfin(q-\alpha,-\alpha,\delta,\gamma,-\varepsilon;1-\midpoint{1}),\\ 
 [f_2,f'_2,r_2,N_2]=\psHeunCsfin(q-\alpha,-\alpha,\delta,\gamma,-\varepsilon;1-\midpoint{1}). 
\end{gather*}
Then we solve the linear system
\[
\left(\begin{matrix} f_1 & f_2\\
 -f'_1 & -f'_2
\end{matrix}\right)
\left(\begin{matrix}
C_1\\C_2
\end{matrix}\right)=
\left(\begin{matrix}
f_0\\f'_0
\end{matrix}\right).
\]
It may also be reasonable to keep the computed values $C_1=C_1(q,\alpha,\gamma,\delta,\varepsilon)$,
$C_2=C_2(q,\alpha,\gamma,\delta,\varepsilon)$ in computer memory.

On finding $C_1$, $C_2$ (by computation or in the computer memory), we define the function
\begin{equation*}
 \psHeunClnearone: z\mapsto [f,f',r,N],
\end{equation*}
where 
\begin{gather*}
 f=C_1 f_1 + C_2 f_2,\quad 
 f'=-C_1 f'_1 - C_2 f'_2,\quad
 r = |C_1|r_1 + |C_2|r_2,\quad
 N = N_1 + N_2,\\
 [f_1,f'_1,r_1,N_1]=\psHeunClfin(q-\alpha,-\alpha,\delta,\gamma,-\varepsilon;1-z),\\
 [f_2,f'_2,r_2,N_2]=\psHeunCsfin(q-\alpha,-\alpha,\delta,\gamma,-\varepsilon;1-z).
\end{gather*}

The described scheme can be repeated literally to define $\psHeunCsnearone$ based
on the representation
\begin{equation}
\HeunCs(q,\alpha,\gamma,\delta,\varepsilon;z) = 
C'_1\HeunCl(q-\alpha,-\alpha,\delta,\gamma,-\varepsilon;1-z)+
C'_2\HeunCs(q-\alpha,-\alpha,\delta,\gamma,-\varepsilon;1-z),
\label{eq:near1s}
\end{equation}
where $C'_1$ and $C'_2$ are some coefficients to be found.

It is notable that finding $C_1$, $C_2$ or $C'_1$, $C'_2$ includes computation of
all three terms in \eqref{eq:near1}
or \eqref{eq:near1s} at $z=\midpoint{1}$. So, if the matching constants are not known, the
algorithms $\psHeunClnearone$ and $\psHeunCsnearone$ are preferable over
$\psHeunClfin$ and $\psHeunCsfin$ in a sufficiently small vicinity of $z=1$. In
the code \cite{mycode2}, used in \S\,\ref{sect:num}, the algorithms are applied for $|z-1|<0.05$.


Consider now a vicinity of the point $z=\infty$. Here the situation is more involved because of the nature of the singular point and in view of the choice of the branch cuts. If the definition of single-valued confluent Heun function demands both branch cuts $\BCzero$,
$\BCone$, then they split the vicinity of infinity
$\{z:|z|>1\}$ into two sectors
$\SectInf^{(\pm)}=\{z:|z|>1,\pm\Im z>0\}$, and coefficients connecting
the function $\HeunCl(z)$ or $\HeunCs(z)$ with two local solutions at infinity are found for
each of the sectors separately.

Assume further in this section that $\varepsilon\neq0$.  For $z\in\SectInf^{\pm}$, we write
\begin{equation}
\HeunC_{\!A,\infty}(q,\alpha,\gamma,\delta,\varepsilon;z)=
E^{\pm}_1 \HeunCl(q,\alpha,\gamma,\delta,\varepsilon;z) + 
E^{\pm}_2 \HeunCs(q,\alpha,\gamma,\delta,\varepsilon;z),
\label{eq:nearinf}
\end{equation}
where $\HeunC_{\!A,\infty}(q,\alpha,\gamma,\delta,\varepsilon;z)$ is the function defined by \eqref{eq:HeunAinfseries} and $E^{\pm}_1$, $E^{\pm}_2$ are some constants. 
It is important that the function in the left-hand side of \eqref{eq:nearinf} does not contain exponential factor while each of the functions in the right-hand side does generally contain (via contribution of $\HeunC_{\!B,\infty}(z)$; see \eqref{eq:expoutinf}). So it is reasonable to choose the matching point to be close to zero and continue $\HeunC_{\!A,\infty}(z)$ from far-field to this point (not vice versa). 
Of importance is also the choice of direction along which the connection of far-field and matching points is realized. We note that Wronskian of $\HeunCl(z)$ and $\HeunCs(z)$, up to a constant factor, is equal to $z^{-\gamma} (1-z)^{-\delta}\E{-\varepsilon z}$ (Liouville--Ostrogradski formula). In view of the exponent, the matrix which arises when finding $E^{\pm}_1$, $E^{\pm}_2$ via matching at a point is usually better conditioned when the point belongs to the so-called anti-Stokes line $\ii \varepsilon^{-1}t$, $t\in(-\infty,\infty)$.

Hence, in the numerical code \cite{mycode2} used in \S~\ref{sect:num}, we choose the following matching point
\[
\midpoint{\pm} = \tfrac{5}{4}\E{\ii\theta^\pm} \quad\mbox{for}\ \
z\in\SectInf^{\pm},
\]
where $\theta^{\pm}=\arg\bigl(\ii\varepsilon^{-1}\bigr)$ if $\ii\varepsilon^{-1}\in\SectInf^{\pm}$ or $\theta^{\pm}=\arg\bigl(-\ii\varepsilon^{-1}\bigr)$ otherwise.
By using the algorithm described in \S~\ref{sect:basic_algorithm}, we continue $\HeunC_{\!A,\infty}(z)$  to the point $\midpoint{\pm}$ starting from $\HeunC_{\!A,\infty}(z)$, $\HeunC'_{\!A,\infty}(z)$ computed with \eqref{eq:HeunAinfseries} at $z^\pm_\infty=R|\varepsilon|^{-1}\E{\ii\theta^\pm}$. The value of $R$ (``far-field radius'') is defined in \cite{mycode2} by the condition that the minimal term in the asymptotic series $\sum_{n=0}^{\infty} n!/R^n$ should be smaller than the machine epsilon $\epsilon$. So we hope that optimal truncation in \eqref{eq:HeunAinfseries} (at series' least term; see, e.g.\ \cite{BerryHowls2015}) for  $z=z_\infty^\pm$ would lead to accuracy of order $\epsilon$.

Thus, we find
\begin{gather*}
 [f^\pm_0,f^\pm_0\vphantom{f}',r^\pm_0,N^\pm_0]=
\psHeunC_{\!\bigl(z^\pm_\infty,\HeunCtxst_{\!A,\infty}(z^\pm_\infty), \HeunCtxst'_{\!A,\infty}(z^\pm_\infty)\bigr)}(\midpoint{\pm}),\\
 [f^\pm_1,f^\pm_1\vphantom{f}',r^\pm_1,N^\pm_1]=\psHeunClfin(\midpoint{\pm}),\qquad
 [f^\pm_2,f^\pm_2\vphantom{f}',r^\pm_2,N^\pm_2]=\psHeunCsfin(\midpoint{\pm}),
\end{gather*}
and the matching coefficients are defined as solution to the linear system
\[
\left(\begin{matrix} f^\pm_1 & f^\pm_2\\
 -f^\pm_1\vphantom{f}' & -f^\pm_2\vphantom{f}'
\end{matrix}\right)
\left(\begin{matrix}
E^\pm_1\\E^\pm_2
\end{matrix}\right)=
\left(\begin{matrix}
f^\pm_0\\f^\pm_0\vphantom{f}'
\end{matrix}\right).
\]

Now, analogously to \eqref{eq:nearinf}, we can write 
\begin{equation}
\HeunC_{\!B,\infty}(q,\alpha,\gamma,\delta,\varepsilon;z)=
D^{\pm}_1 \HeunCl(q,\alpha,\gamma,\delta,\varepsilon;z) + 
D^{\pm}_2 \HeunCs(q,\alpha,\gamma,\delta,\varepsilon;z),
\label{eq:nearinfB}
\end{equation}
where $\HeunC_{\!B,\infty}(q,\alpha,\gamma,\delta,\varepsilon;z)$ is the function defined by \eqref{eq:expoutinf} and $D^{\pm}_1$, $D^{\pm}_2$ are some constants. However, unlike the previous case all terms of \eqref{eq:nearinfB} may have exponential factor $\E{-\varepsilon z}$. So, it is useful to transform equality \eqref{eq:nearinfB} using relationship \eqref{eq:expout}. We write
\begin{align*}\notag
\HeunC_{\!A,\infty}(q-\varepsilon\gamma,\alpha-\varepsilon(\gamma+\delta),\gamma,\delta,-\varepsilon;z)={}&
D^{\pm}_1 \HeunCl(q-\varepsilon\gamma,\alpha-\varepsilon(\gamma+\delta),\gamma,\delta,-\varepsilon;z) \\&{}+ 
D^{\pm}_2 \HeunCs(q-\varepsilon\gamma,\alpha-\varepsilon(\gamma+\delta),\gamma,\delta,-\varepsilon;z),
\label{eq:nearinfBmod}
\end{align*}
and the procedure described above can be applied literally.

Introduce now the matrix $\bm{d}^\pm=\bigl[d^\pm_{ij}\bigr]_{i,j=1,2}$:
\[
\bm{d}^\pm=\left(\begin{matrix}
E^\pm_1 & E^\pm_2\\
D^\pm_1 & D^\pm_2\\
\end{matrix}\right)^{-1}.
\]
Then, for large $|\varepsilon z|$ ($|\varepsilon z|>R$ in code \cite{mycode2} used in \S~\ref{sect:num}), we compute
\begin{gather*}
 \HeunCl(q,\alpha,\gamma,\delta,\varepsilon;z) = 
 d^\pm_{11}\HeunC_{\!A,\infty}(q,\alpha,\gamma,\delta,\varepsilon;z)+
 d^\pm_{12}\HeunC_{\!B,\infty}(q,\alpha,\gamma,\delta,\varepsilon;z),
 \label{eq:nearinftyl}
\\
 \HeunCs(q,\alpha,\gamma,\delta,\varepsilon;z) = 
 d^\pm_{21}\HeunC_{\!A,\infty}(q,\alpha,\gamma,\delta,\varepsilon;z)+
 d^\pm_{22}\HeunC_{\!B,\infty}(q,\alpha,\gamma,\delta,\varepsilon;z),
 \label{eq:nearinftys}
\end{gather*}
when $z\in\SectInf^\pm$. We denote by $\psHeunClnearinf$ and $\psHeunCsnearinf$ the described algorithms for finding  
$\HeunCl(z)$ and $\HeunCs(z)$ for large $|\varepsilon z|$.
 
\section{Numerical results}
\label{sect:num}

In this section we present results of numerical evaluation of the functions
$\HeunCl(z)$ and $\HeunCs(z)$. For these tests we use both the basic algorithms $\psHeunClfin$, $\psHeunCsfin$ (see \S~\ref{sect:basic_algorithm}) and the algorithms with improvements described in the previous section ($\psHeunClnearone$, $\psHeunCsnearone$, $\psHeunClnearinf$, $\psHeunCsnearinf$).

Calculations are performed with the code \cite{mycode2} in the numerical computing environment GNU Octave and double precision (64-bit)
arithmetics (the machine epsilon $\epsilon$ is about $2.22\cdot10^{-16}$).

It is rather straightforward to check the following special forms of the confluent Heun functions
\begin{equation*}
\HeunCl\bigl(\tfrac14,0,\tfrac12,\tfrac12,0;z\bigr)=\sqrt{1-z}:=h_1(z),\qquad 
\HeunCs\bigl(\tfrac14,0,\tfrac12,\tfrac12,0;z\bigr)=\sqrt{z}:=h_2(z),
\end{equation*}
and, so, we define
\begin{gather*}
\Delta_1(z):=\HeunCl\bigl(\tfrac14,0,\tfrac12,\tfrac12,0;z\bigr)-h_1(z),\\
\Delta_2(z):=\HeunCs\bigl(\tfrac14,0,\tfrac12,\tfrac12,0;z\bigr)-h_2(z).
\end{gather*}
Further we will check numerically the identities $\Delta_n(z)=0$, $n=1,2,\ldots,9$, where, analogously, we introduce
\begin{gather*}
\Delta_3(z):=\HeunCl(6,0,1,1,0;z)-h_3(z),\qquad h_3(z):=6 z^2-6 z+1,\\
\Delta_4(z):=\HeunCs(6,0,1,1,0;z)-h_4(z),\qquad h_4(z):=\bigl(6 z^2-6 z+1\bigr)\bigl(\log(z)-\log(1-z)-3\bigr)-6 z+3,\\
\Delta_5(z):=\HeunCl\bigl(-\tfrac14,0,\tfrac12,\tfrac12,0;z\bigr)-h_5(z),\qquad h_5(z):=\cos\log\bigl(\sqrt{1-z}+\ii\sqrt{z}\bigr),\\
\Delta_6(z):=\HeunCs\bigl(-\tfrac14,0,\tfrac12,\tfrac12,0;z\bigr)-h_6(z),\qquad h_6(z):=-\ii\sin\log\bigl(\sqrt{1-z}+\ii\sqrt{z}\bigr),\\
\Delta_7(z):=\HeunCl\bigl(\tfrac34,\tfrac32,\tfrac12,\tfrac12,1;z\bigr)-h_7(z),\qquad h_7(z):=\E{-z}\sqrt{1-z},\\
\Delta_8(z):=\HeunCs\bigl(\tfrac54,\tfrac32,\tfrac12,\tfrac12,1;z\bigr)-h_8(z),\qquad h_8(z):=\E{-z}\sqrt{z},
\end{gather*}
and
\begin{equation*}
\Delta_9(z):=\HeunCl(-2,0,-1,0,1;z)+\tfrac32\HeunCs(-2,0,-1,0,1;z)-h_9(z),\qquad h_9(z):=\E{-z}(1-z).
\label{eq:Delta9}
\end{equation*}
The coefficient $3/2$ in the definition of $\Delta_9(z)$ can be easily found by comparing expansions at $z=0$ for the three terms in the right-hand side. It should also be marked that $\Delta_4(z)$ and $\Delta_9(z)$ relate to the special cases: when $\gamma=1$ and $\gamma$ is a non-positive integer. 

\begin{figure}[t!]
\centering
 \SetLabels
 \L (-0.0425*0.96) (a)\\
 \L (0.521*0.96) (b)\\
 \L (0*0.86) $\Im z$\\
 \L (0.27*-0.06) $\Re z$\\
 \L (0.5675*0.86) $\Im z$\\
 \L (0.843*-0.06) $\Re z$\\
 \endSetLabels
 \leavevmode\AffixLabels{\mbox{\kern5mm}\includegraphics[height=46mm]{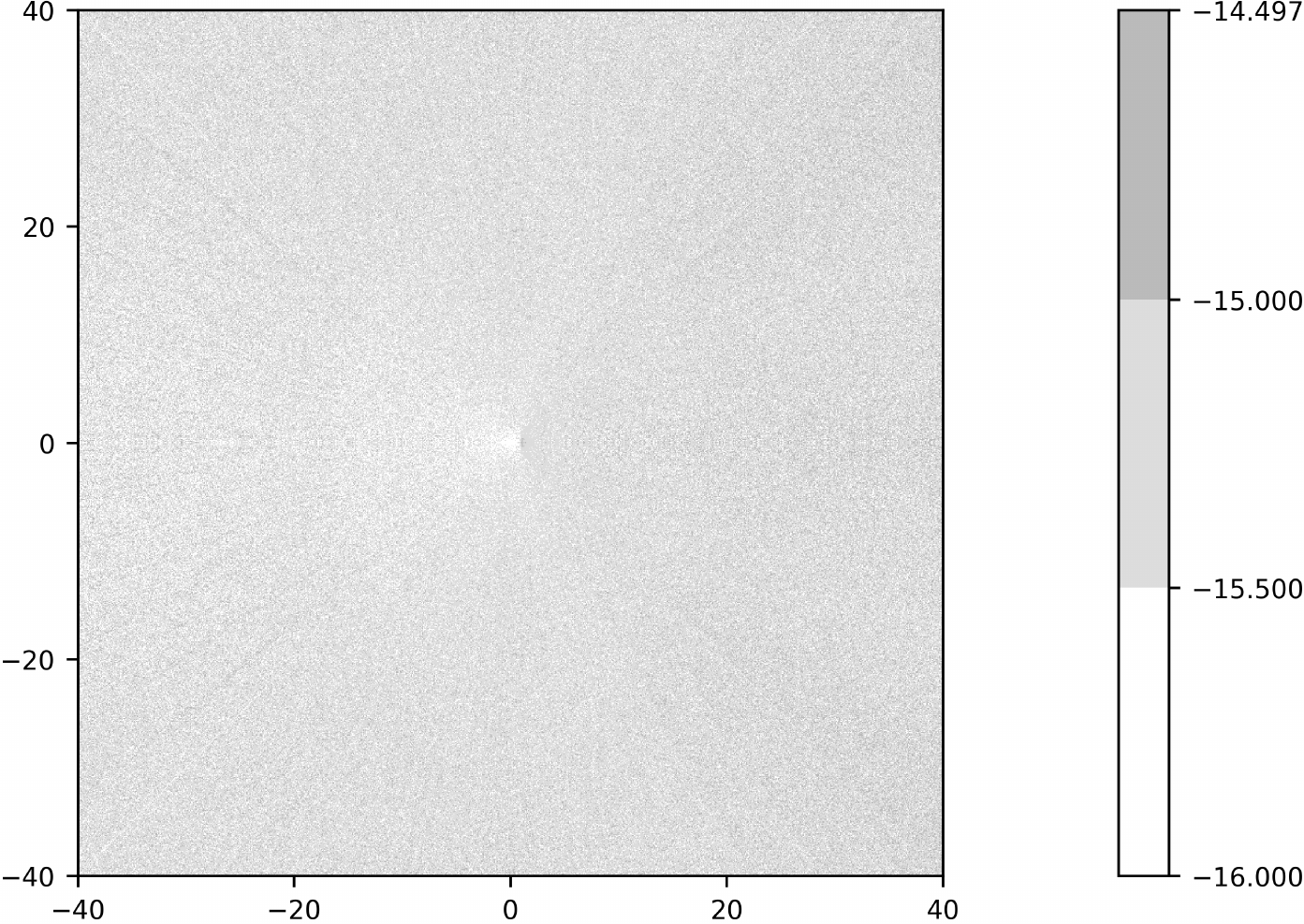}\kern25mm
 \includegraphics[height=46mm]{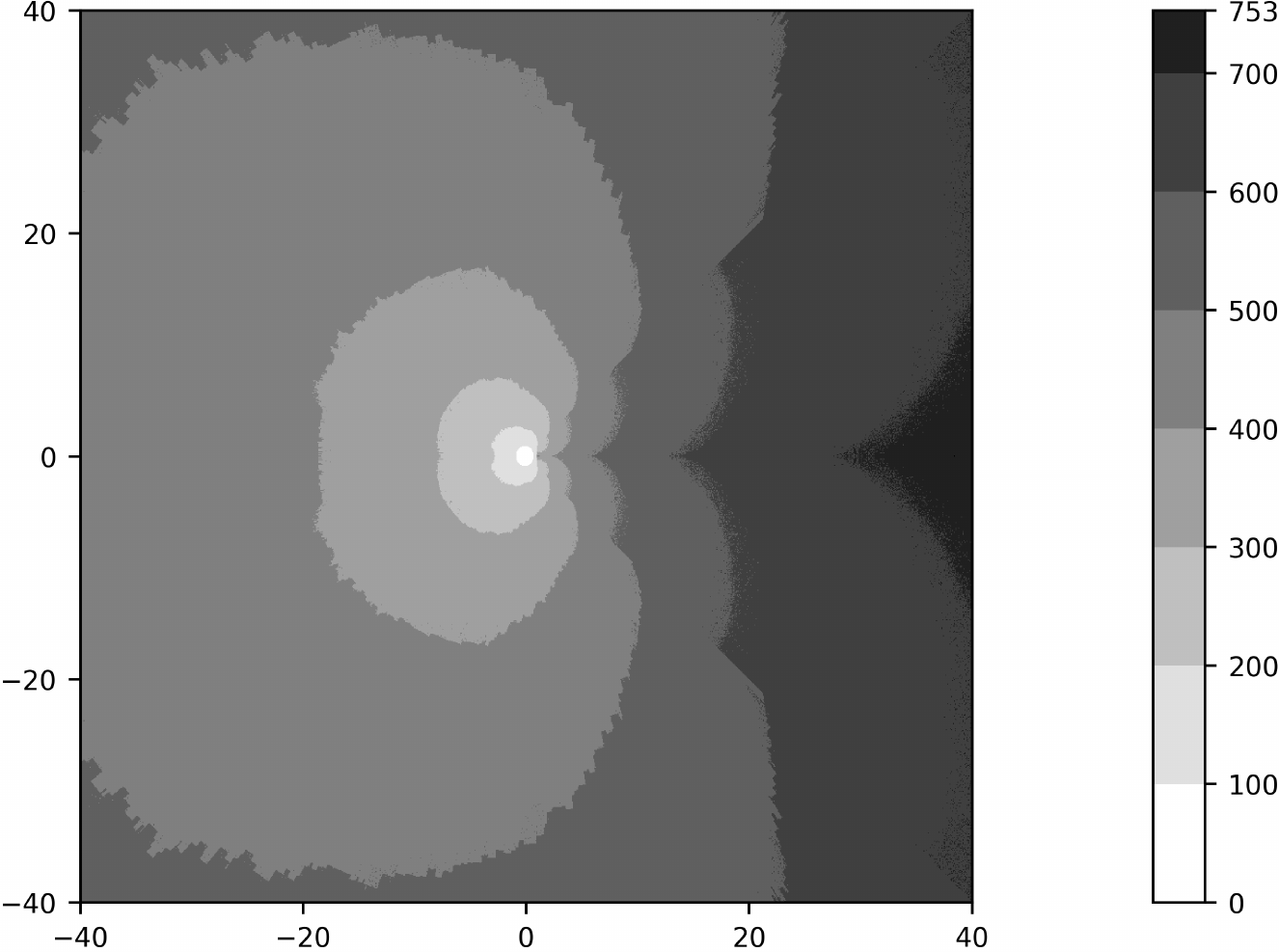}}
 \vspace{2mm}
 \caption{Values of $\max\bigl\{-16,\log_{10}\Lambda_1(z)\bigr\}$ (a) and $N_1(z)$ (b).}
 \label{fig:L1_0}
\centering\vspace{4mm}
 \SetLabels
 \L (-0.0425*0.96) (a)\\
 \L (0.521*0.96) (b)\\
 \L (0*0.86) $\Im z$\\
 \L (0.27*-0.06) $\Re z$\\
 \L (0.5675*0.86) $\Im z$\\
 \L (0.843*-0.06) $\Re z$\\
 \endSetLabels
 \leavevmode\AffixLabels{\mbox{\kern5mm}\includegraphics[height=46mm]{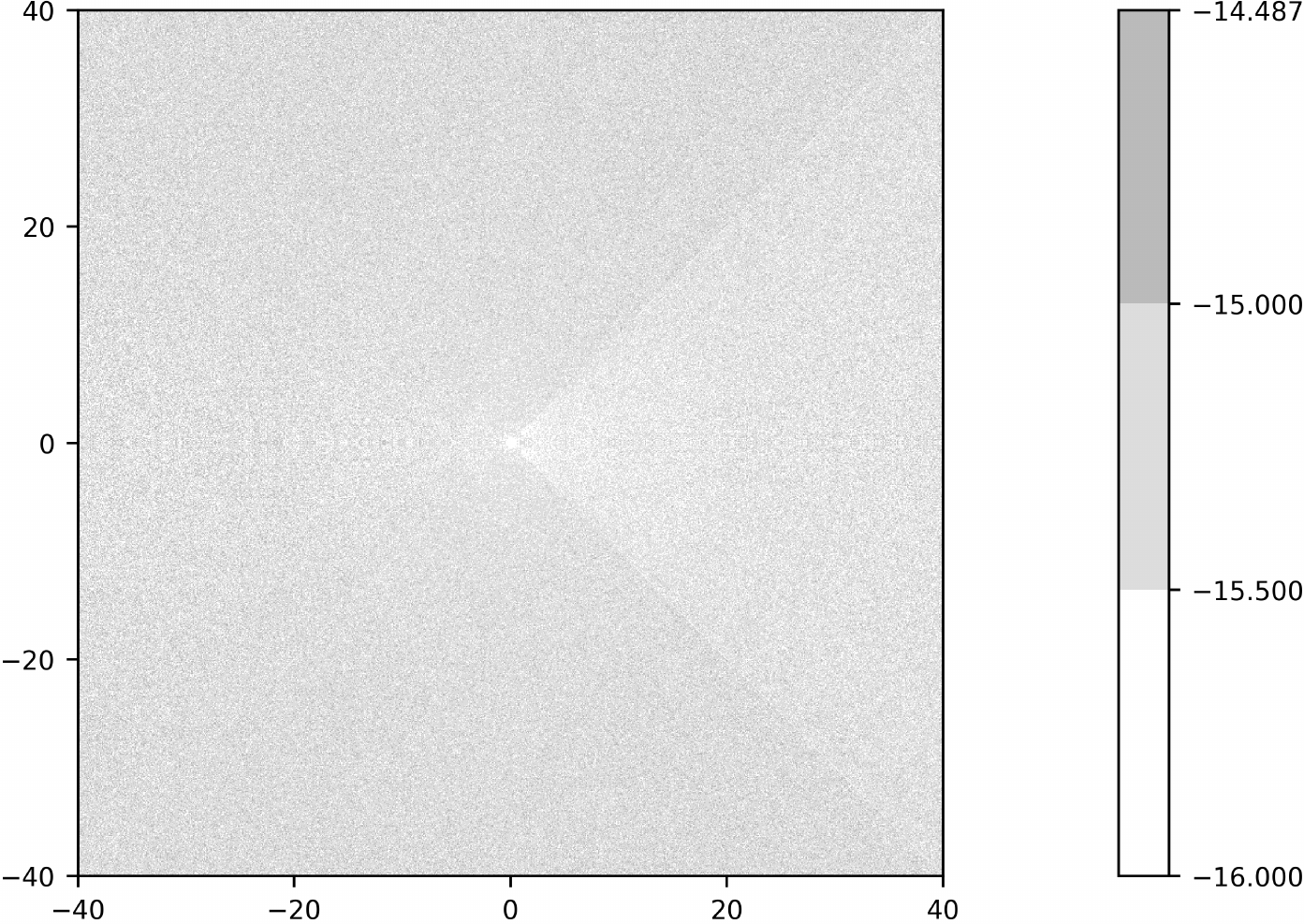}\kern25mm
 \includegraphics[height=46mm]{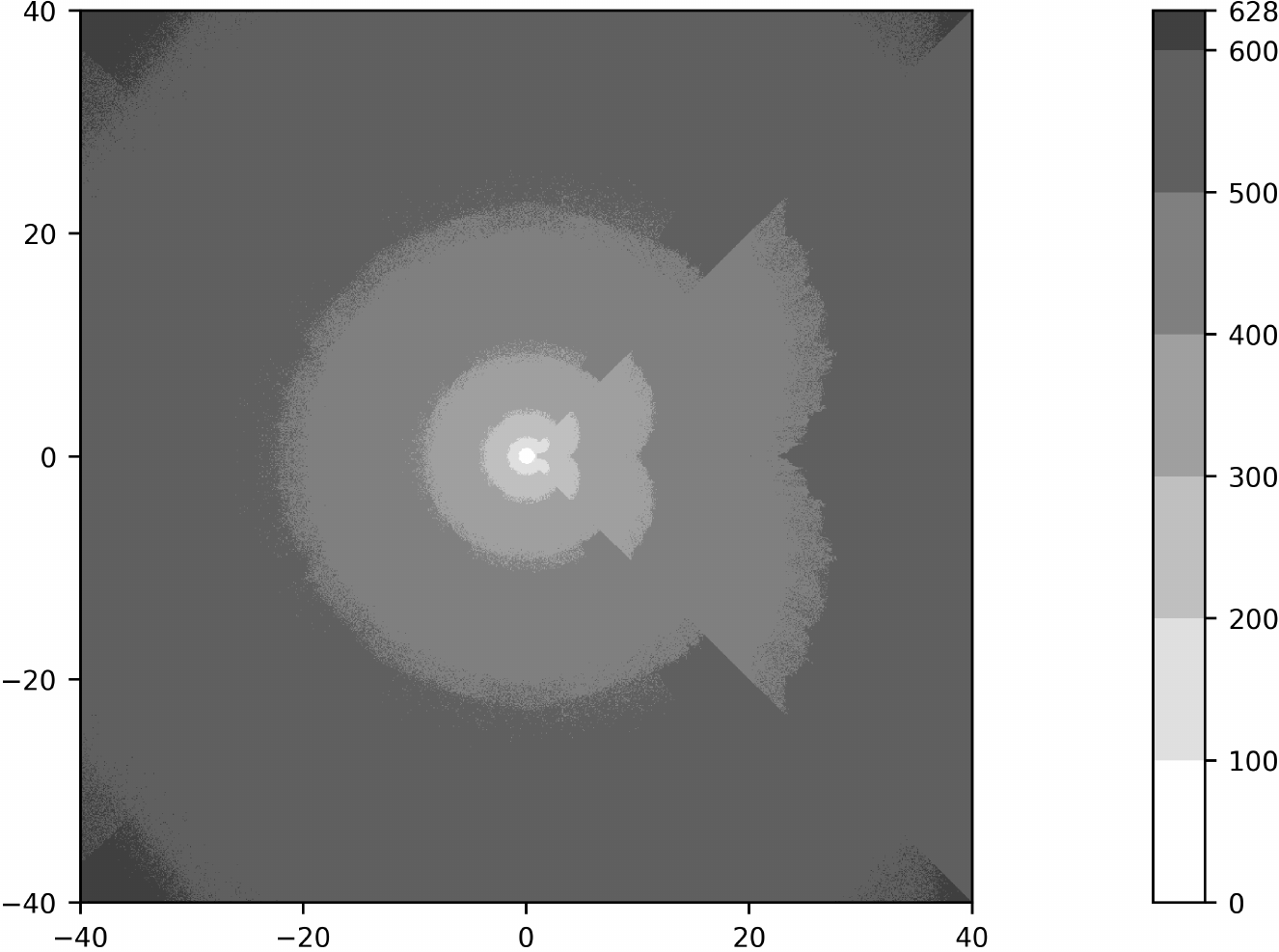}}
 \vspace{2mm}
 \caption{Values of $\max\bigl\{-16,\log_{10}\Lambda_2(z)\bigr\}$ (a) and $N_2(z)$ (b).}
 \label{fig:L2_0}
\end{figure}
\begin{figure}[p!]
\centering
 \SetLabels
 \L (-0.0425*0.96) (a)\\
 \L (0.521*0.96) (b)\\
 \L (0*0.86) $\Im z$\\
 \L (0.27*-0.06) $\Re z$\\
 \L (0.5675*0.86) $\Im z$\\
 \L (0.843*-0.06) $\Re z$\\
 \endSetLabels
 \leavevmode\AffixLabels{\mbox{\kern5mm}\includegraphics[height=46mm]{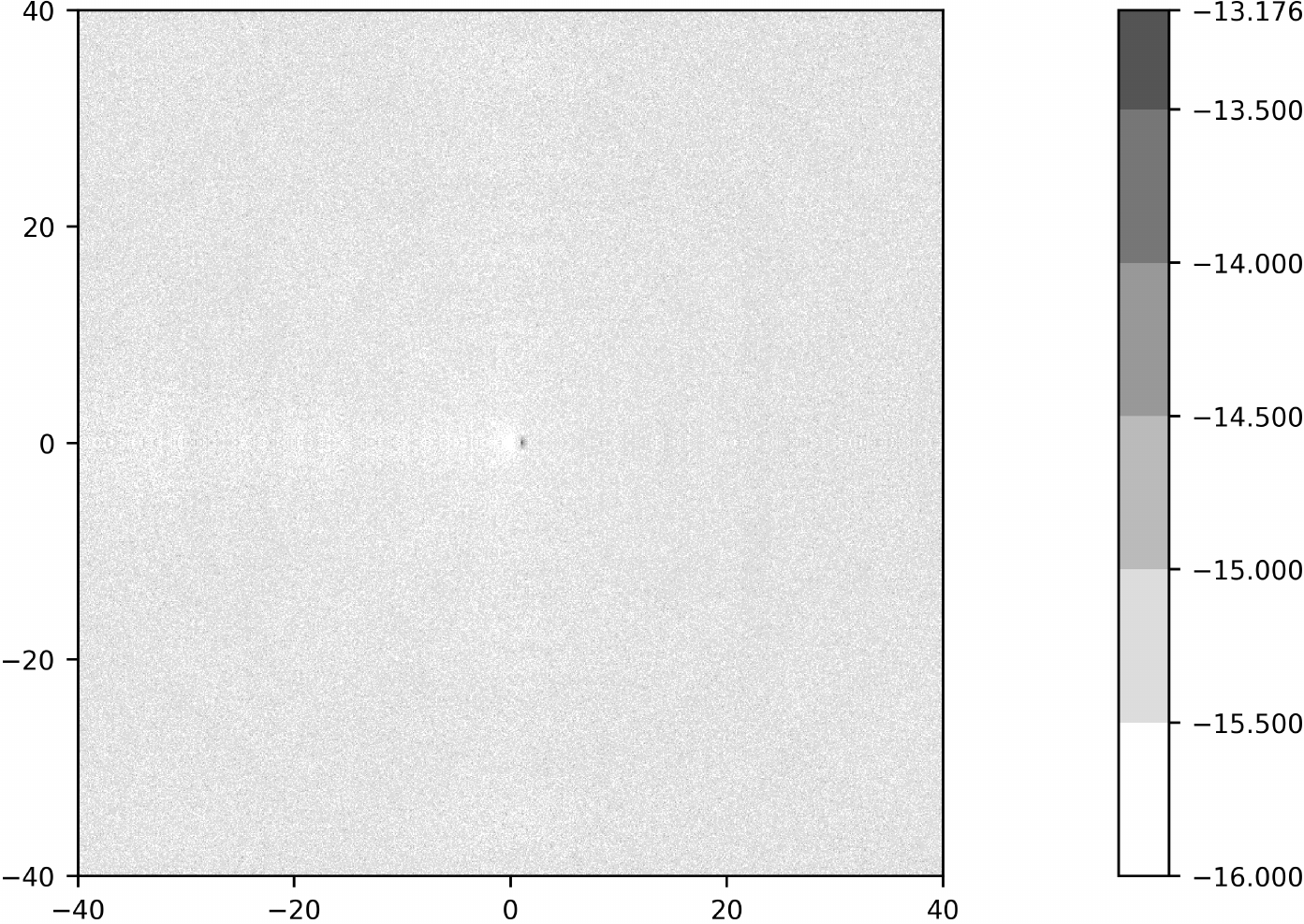}\kern25mm
 \includegraphics[height=46mm]{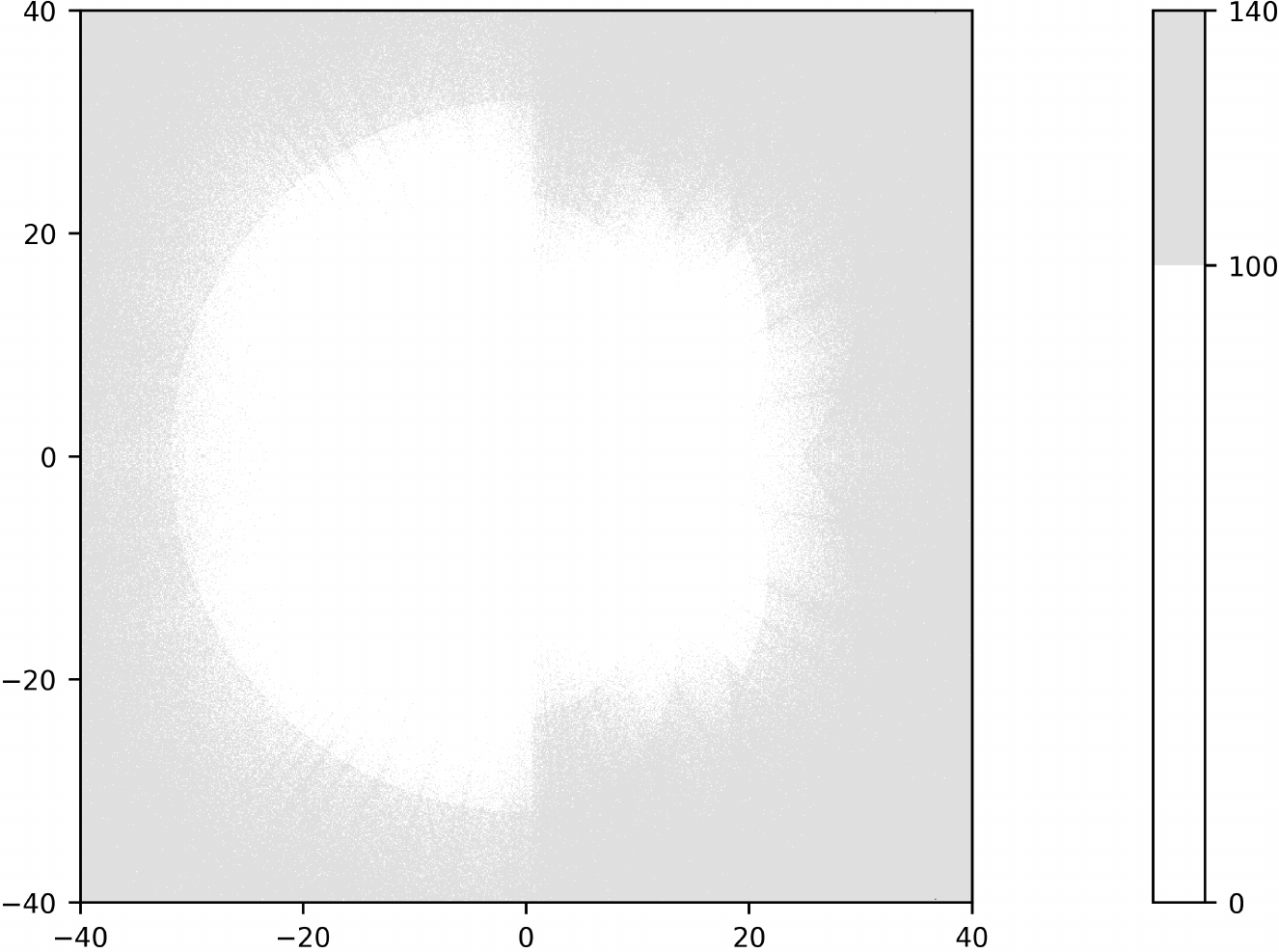}}
 \vspace{2mm}
 \caption{Values of $\max\bigl\{-16,\log_{10}\Lambda_3(z)\bigr\}$ (a) and $N_3(z)$ (b).}
 \label{fig:L3_0}
\centering\vspace{4mm}
 \SetLabels
 \L (-0.0425*0.96) (a)\\
 \L (0.521*0.96) (b)\\
 \L (0*0.86) $\Im z$\\
 \L (0.27*-0.06) $\Re z$\\
 \L (0.5675*0.86) $\Im z$\\
 \L (0.843*-0.06) $\Re z$\\
 \endSetLabels
 \leavevmode\AffixLabels{\mbox{\kern5mm}\includegraphics[height=46mm]{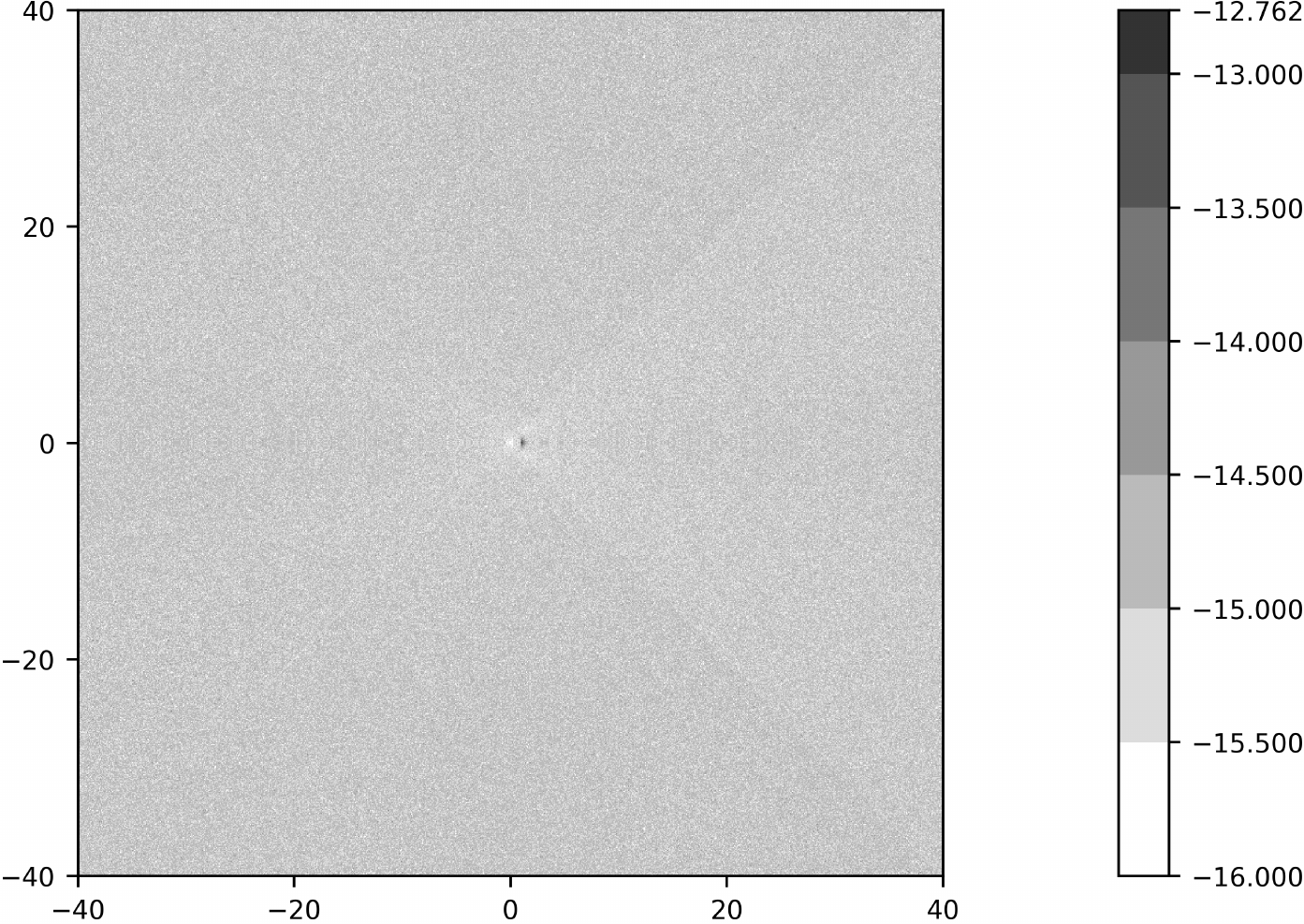}\kern25mm
 \includegraphics[height=46mm]{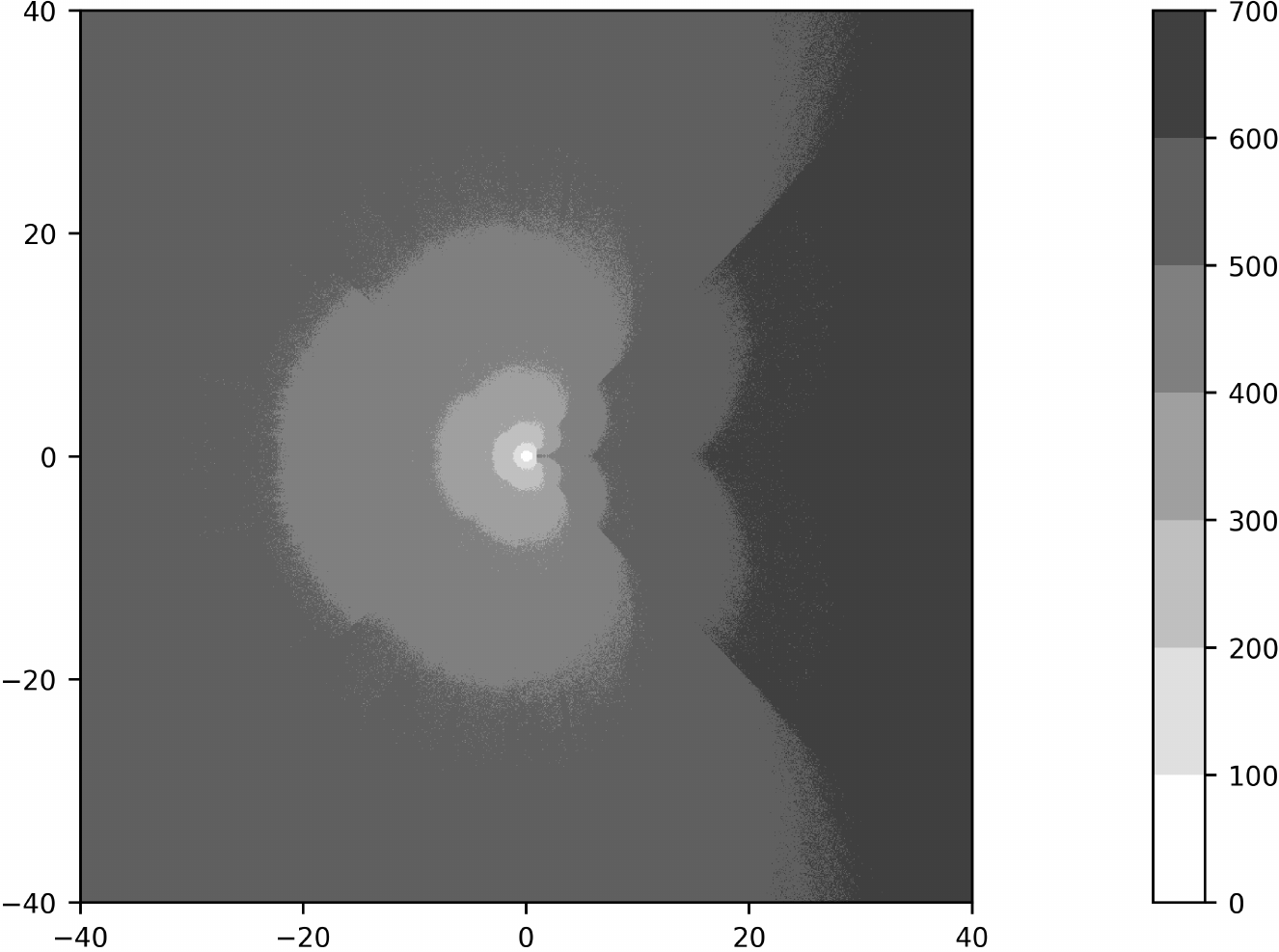}}
 \vspace{2mm}
 \caption{Values of $\max\bigl\{-16,\log_{10}\Lambda_4(z)\bigr\}$ (a) and $N_4(z)$ (b).}
 \label{fig:L4_0}
\centering\vspace{4mm}
 \SetLabels
 \L (-0.0425*0.96) (a)\\
 \L (0.521*0.96) (b)\\
 \L (0*0.86) $\Im z$\\
 \L (0.27*-0.06) $\Re z$\\
 \L (0.5675*0.86) $\Im z$\\
 \L (0.843*-0.06) $\Re z$\\
 \endSetLabels
 \leavevmode\AffixLabels{\mbox{\kern5mm}\includegraphics[height=46mm]{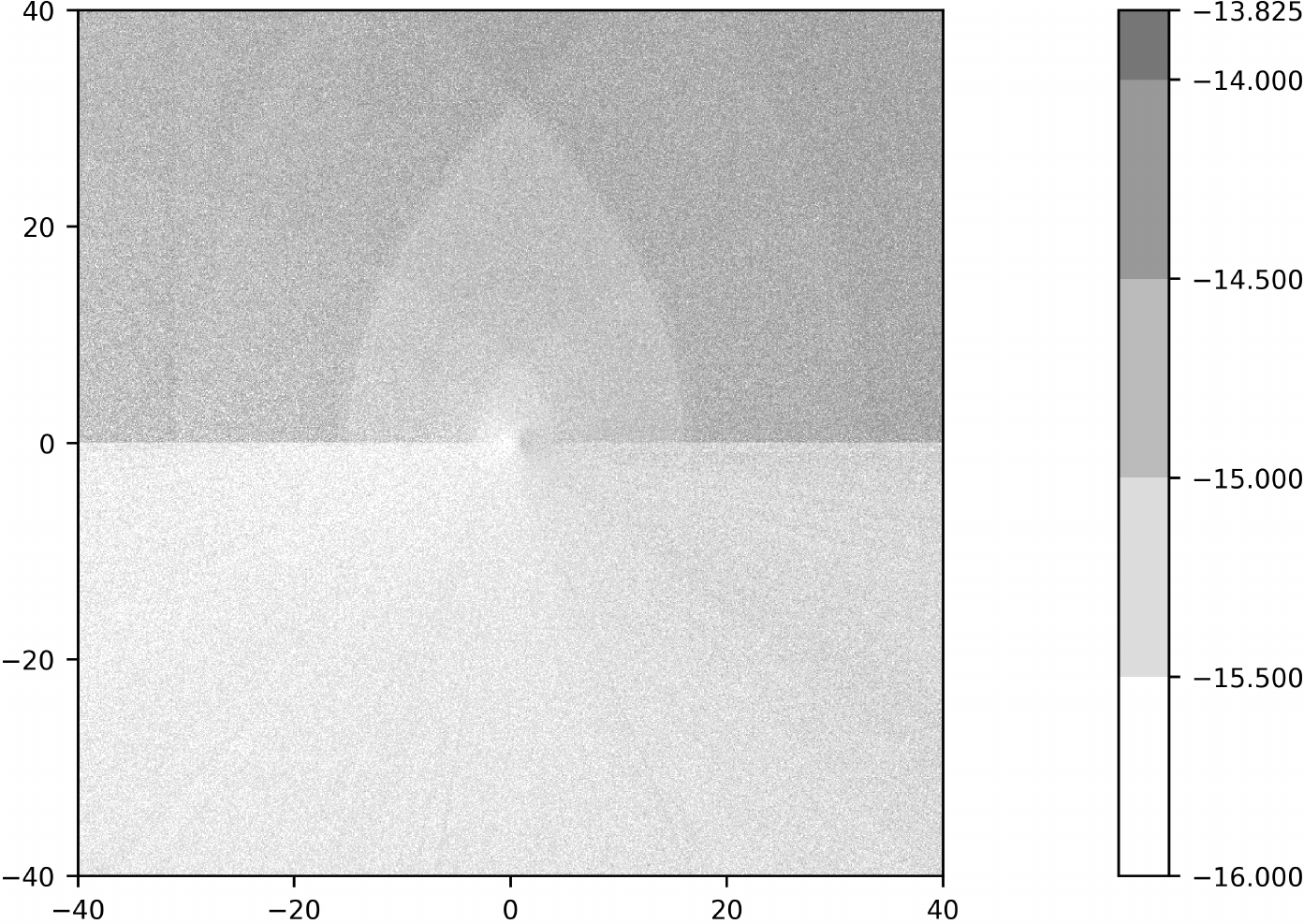}\kern25mm
 \includegraphics[height=46mm]{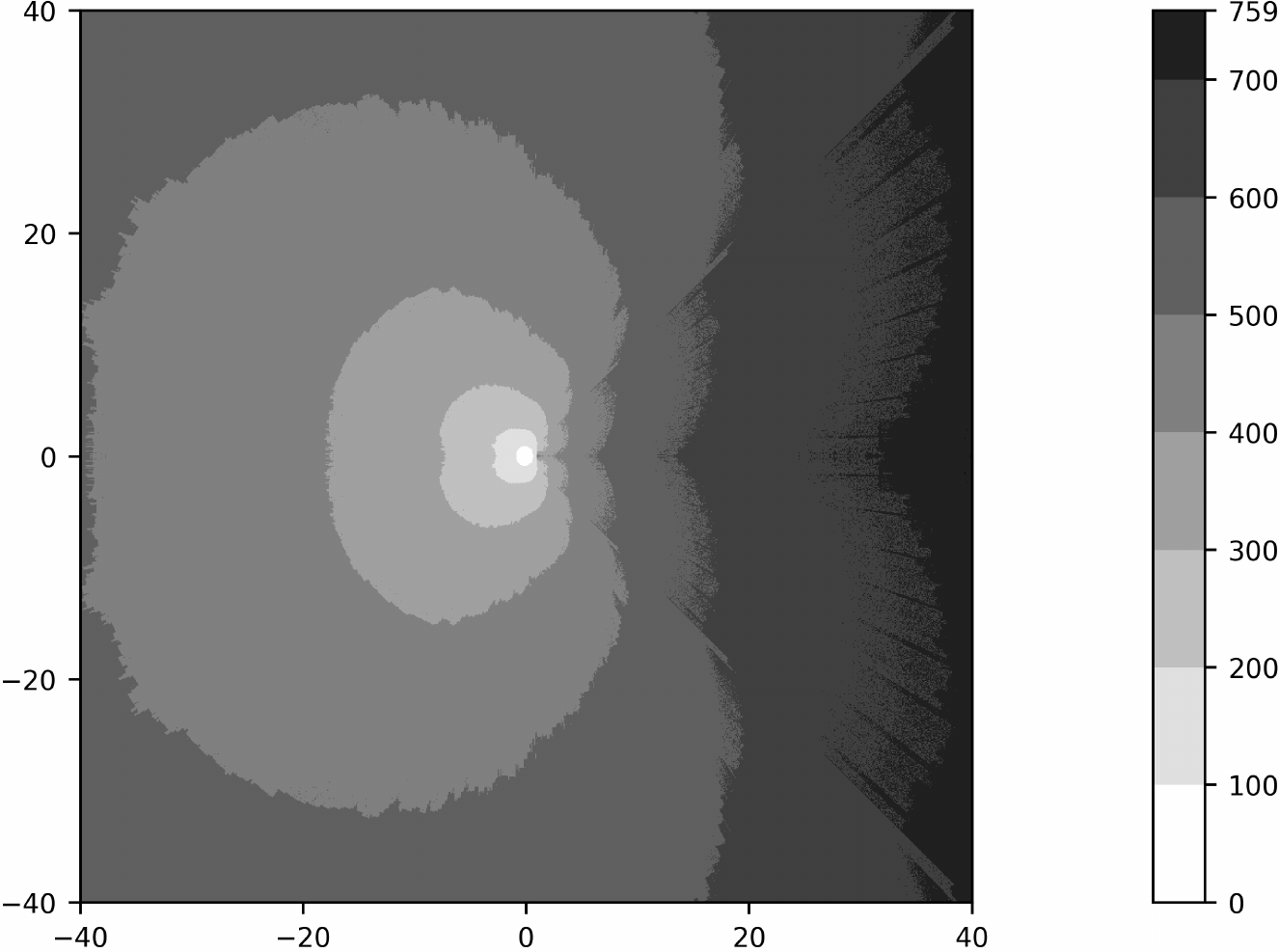}}
 \vspace{2mm}
 \caption{Values of $\max\bigl\{-16,\log_{10}\Lambda_5(z)\bigr\}$ (a) and $N_5(z)$ (b).}
 \label{fig:L5_0}
\centering\vspace{4mm}
 \SetLabels
 \L (-0.0425*0.96) (a)\\
 \L (0.521*0.96) (b)\\
 \L (0*0.86) $\Im z$\\
 \L (0.27*-0.06) $\Re z$\\
 \L (0.5675*0.86) $\Im z$\\
 \L (0.843*-0.06) $\Re z$\\
 \endSetLabels
 \leavevmode\AffixLabels{\mbox{\kern5mm}\includegraphics[height=46mm]{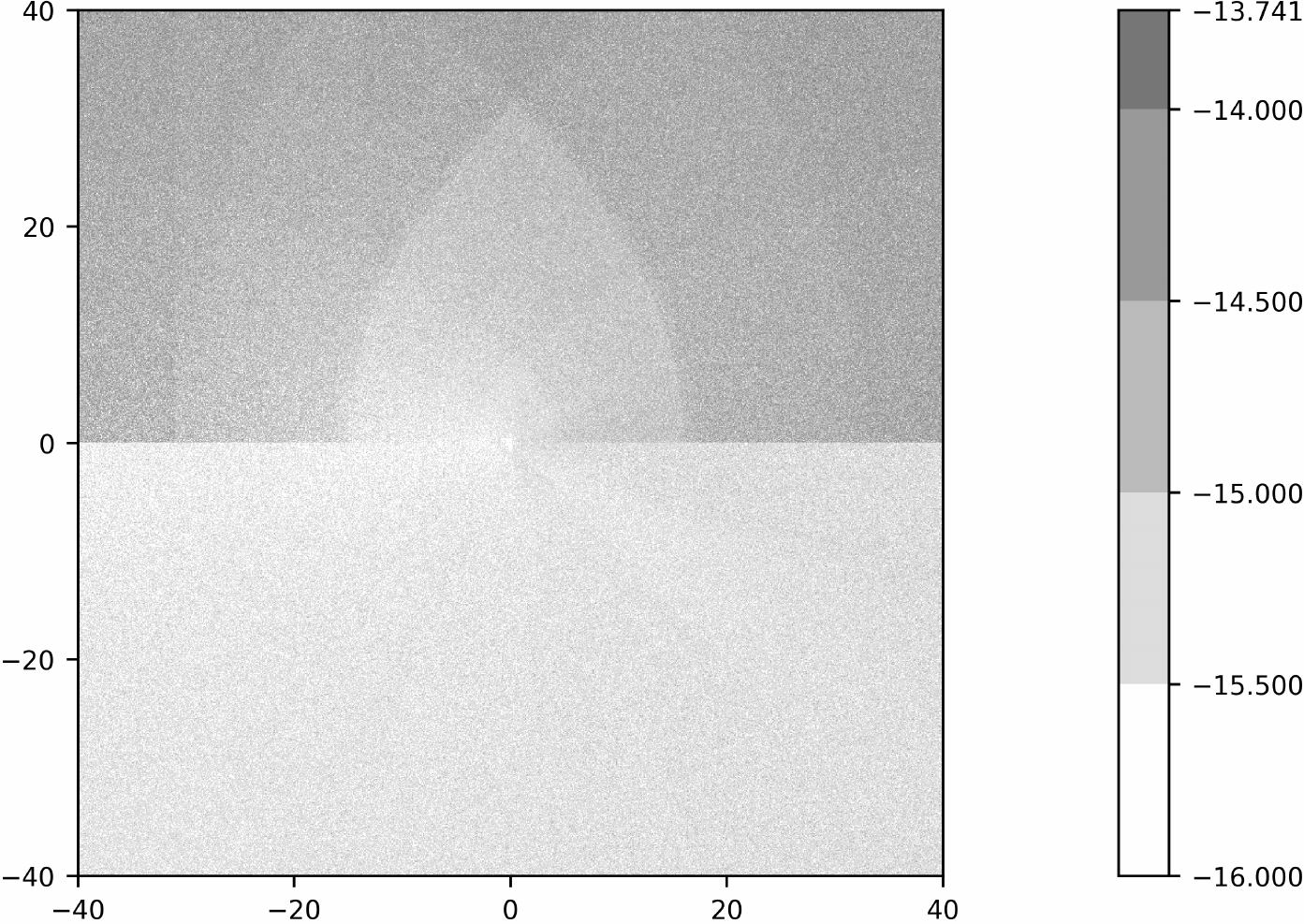}\kern25mm
 \includegraphics[height=46mm]{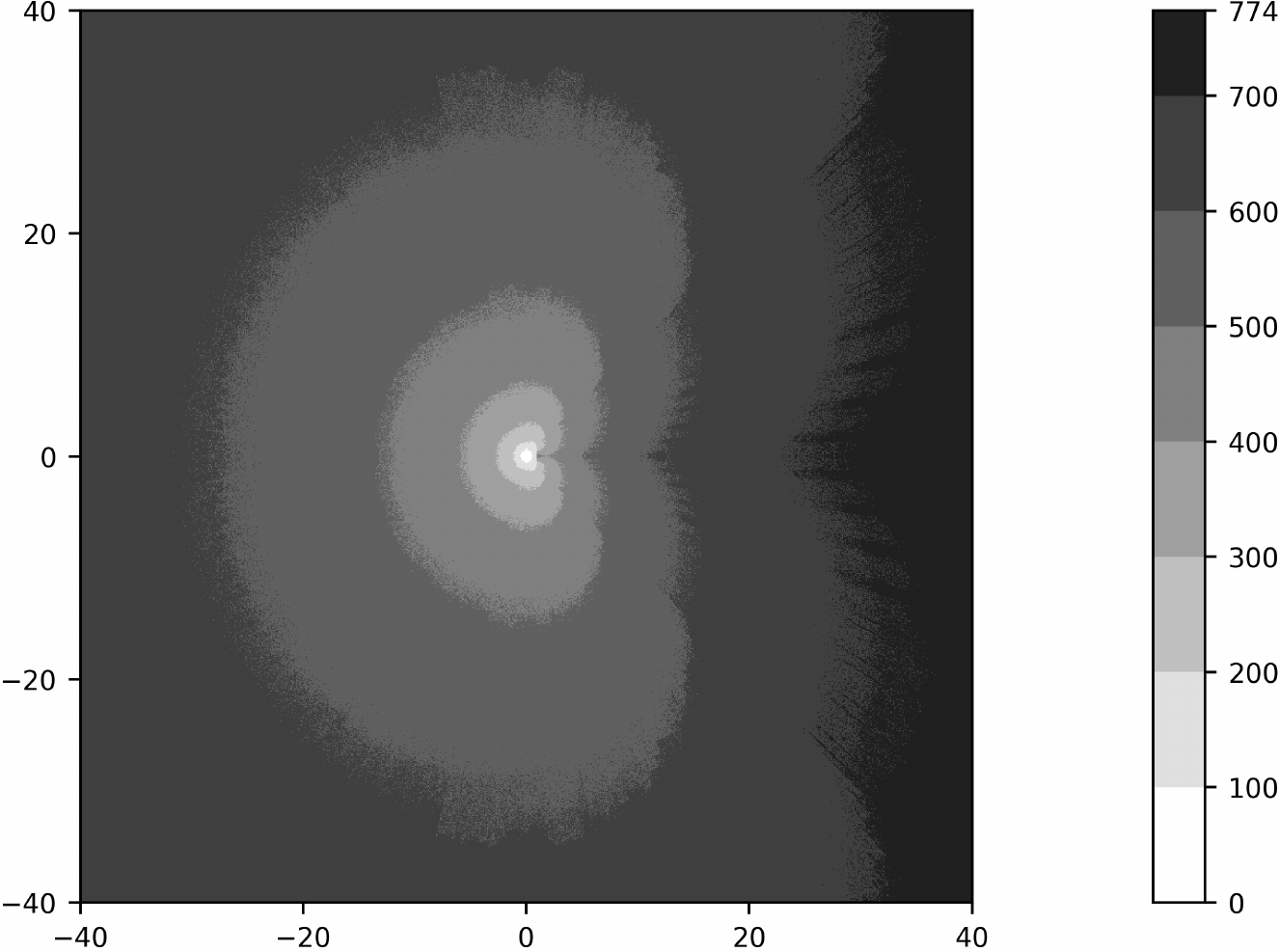}}
 \vspace{2mm}
 \caption{Values of $\max\bigl\{-16,\log_{10}\Lambda_6(z)\bigr\}$ (a) and $N_6(z)$ (b).}
 \label{fig:L6_0}
\end{figure}

\begin{figure}[p!]
\centering
 \SetLabels
 \L (-0.0425*0.96) (a)\\
 \L (0.521*0.96) (b)\\
 \L (0*0.86) $\Im z$\\
 \L (0.27*-0.06) $\Re z$\\
 \L (0.5675*0.86) $\Im z$\\
 \L (0.843*-0.06) $\Re z$\\
 \endSetLabels
 \leavevmode\AffixLabels{\mbox{\kern5mm}\includegraphics[height=46mm]{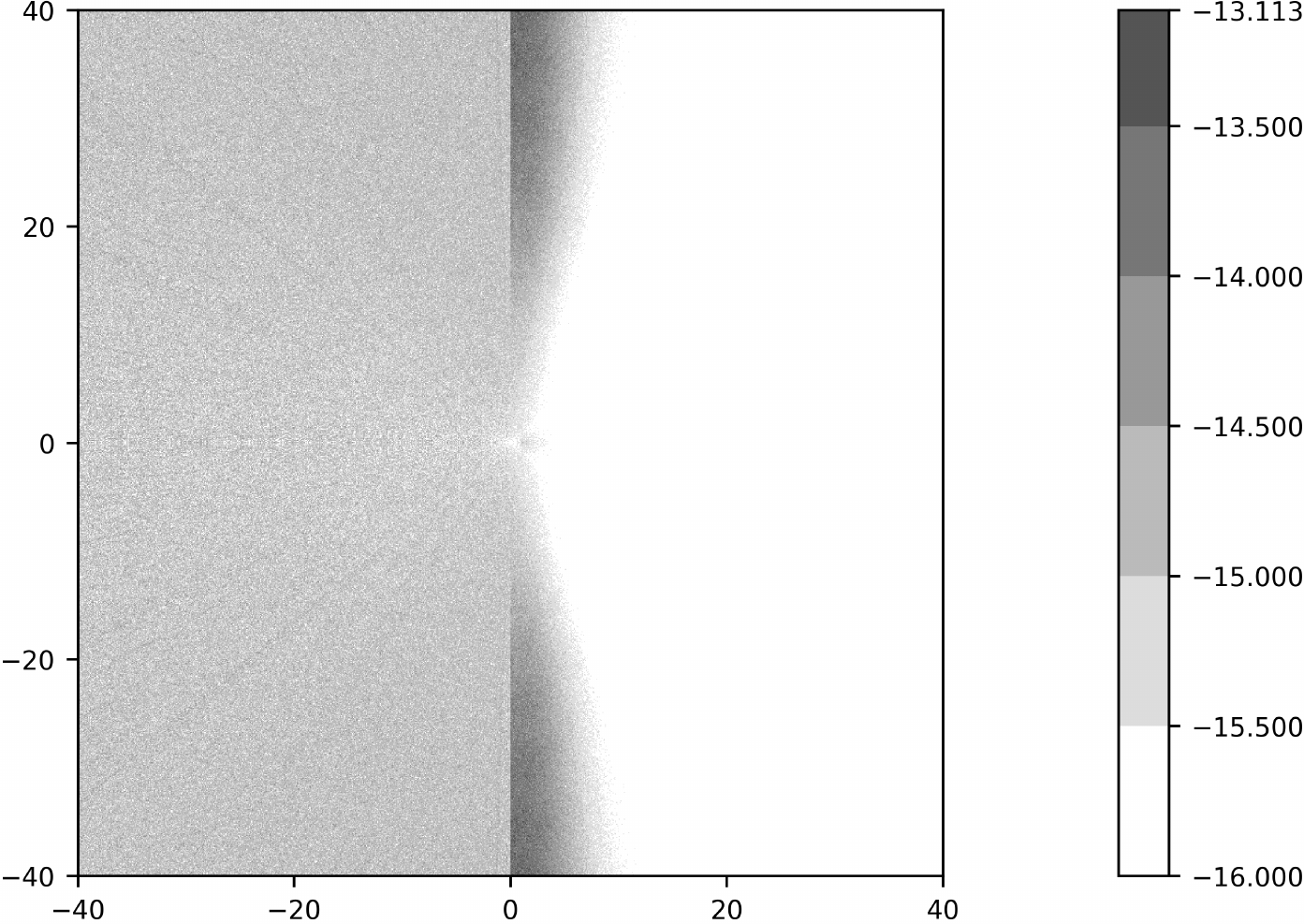}\kern25mm
 \includegraphics[height=46mm]{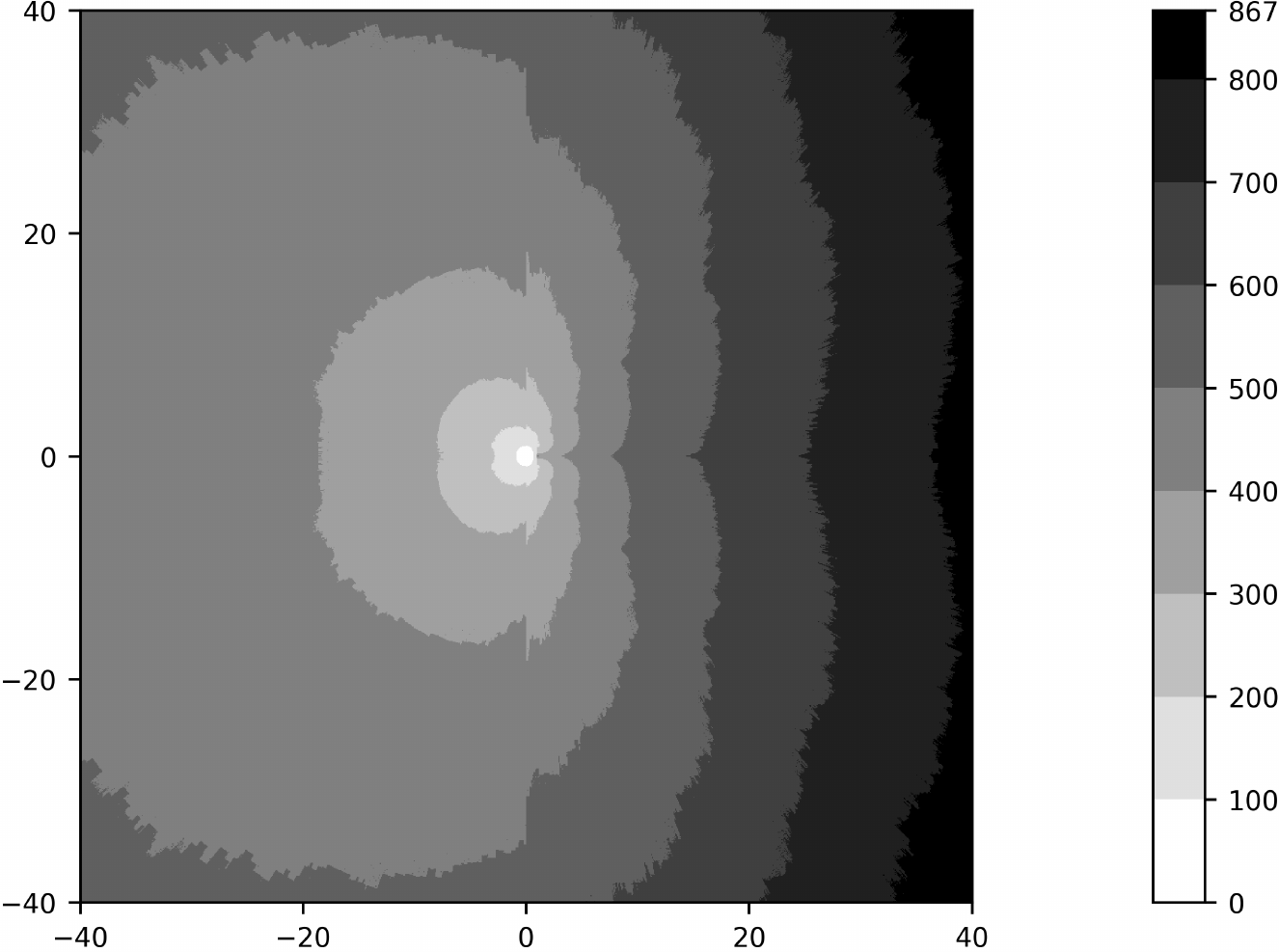}}
 \vspace{2mm}
 \caption{Values of $\max\bigl\{-16,\log_{10}\Lambda_7(z)\bigr\}$ (a) and $N_7(z)$ (b).}
 \label{fig:L7_0}
\centering\vspace{4mm}
 \SetLabels
 \L (-0.0425*0.96) (a)\\
 \L (0.521*0.96) (b)\\
 \L (0*0.86) $\Im z$\\
 \L (0.27*-0.06) $\Re z$\\
 \L (0.5675*0.86) $\Im z$\\
 \L (0.843*-0.06) $\Re z$\\
 \endSetLabels
 \leavevmode\AffixLabels{\mbox{\kern5mm}\includegraphics[height=46mm]{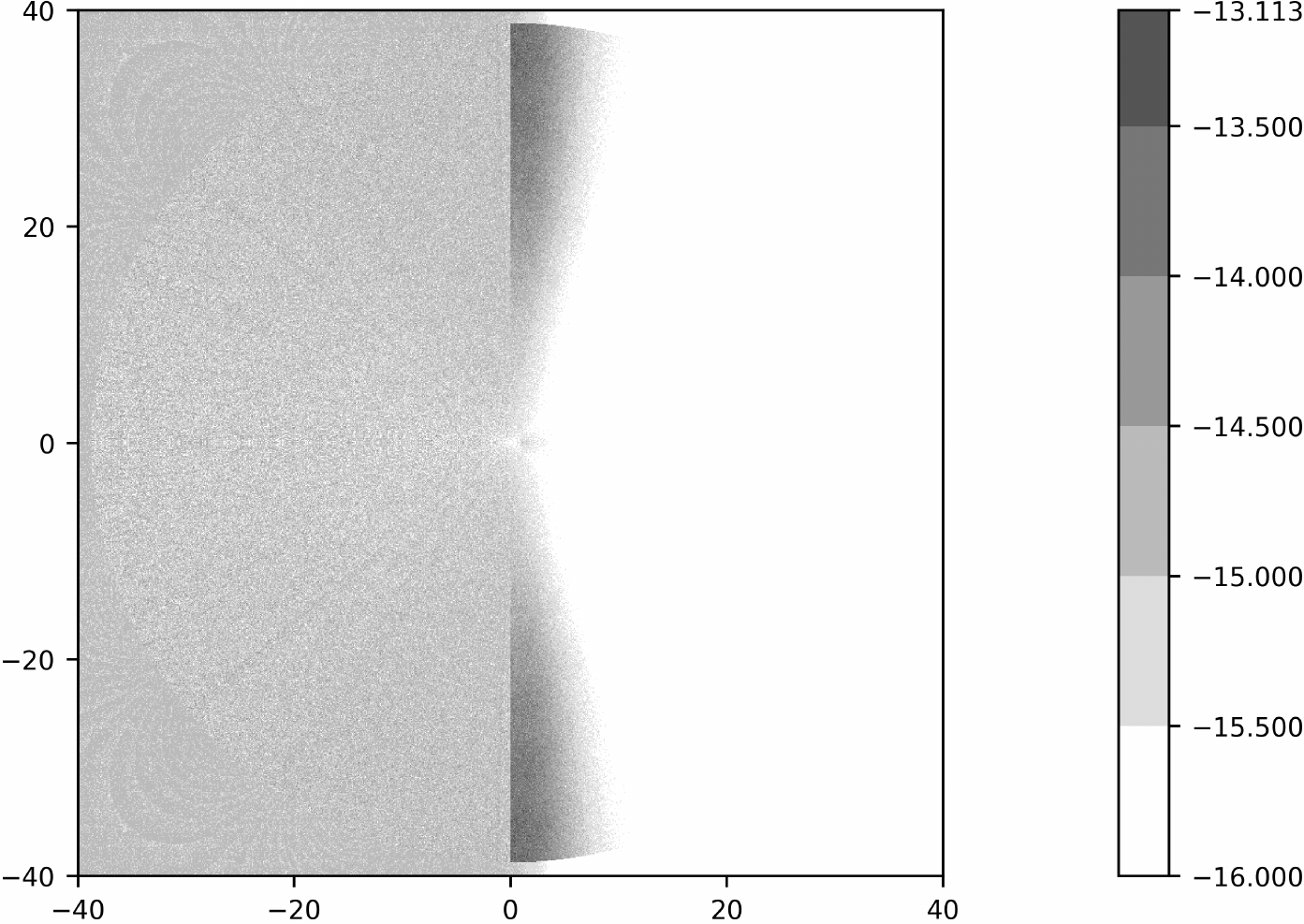}\kern25mm
 \includegraphics[height=46mm]{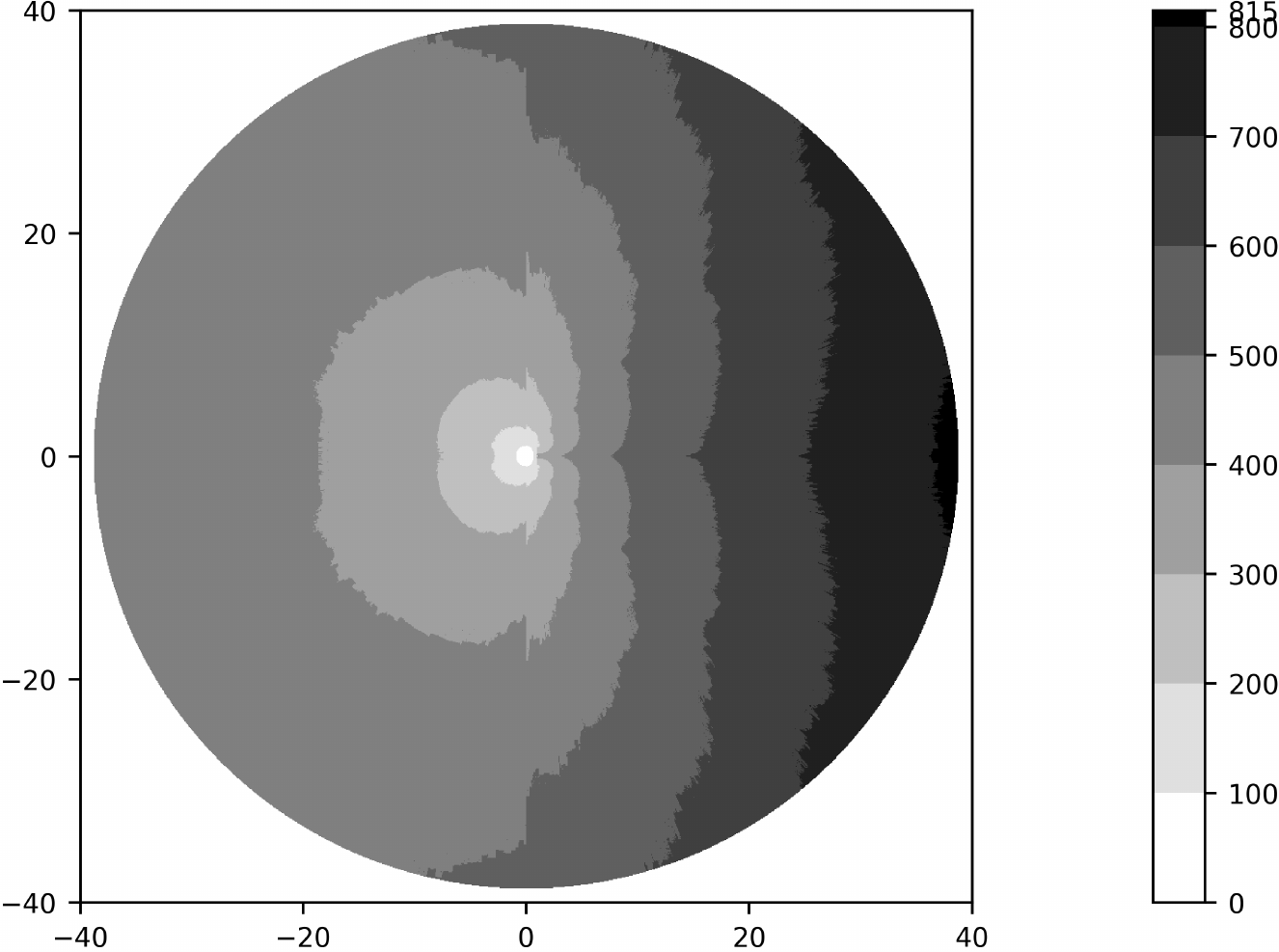}}
 \vspace{2mm}
 \caption{Values of $\max\bigl\{-16,\log_{10}\Lambda_7(z)\bigr\}$ (a) and $N_7(z)$ (b). Improved algorithms $\psHeunClnearone$ and $\psHeunClnearinf$ are used.}
 \label{fig:L7}
\centering\vspace{4mm}
 \SetLabels
 \L (-0.0425*0.96) (a)\\
 \L (0.521*0.96) (b)\\
 \L (0*0.86) $\Im z$\\
 \L (0.27*-0.06) $\Re z$\\
 \L (0.5675*0.86) $\Im z$\\
 \L (0.843*-0.06) $\Re z$\\
 \endSetLabels
 \leavevmode\AffixLabels{\mbox{\kern5mm}\includegraphics[height=46mm]{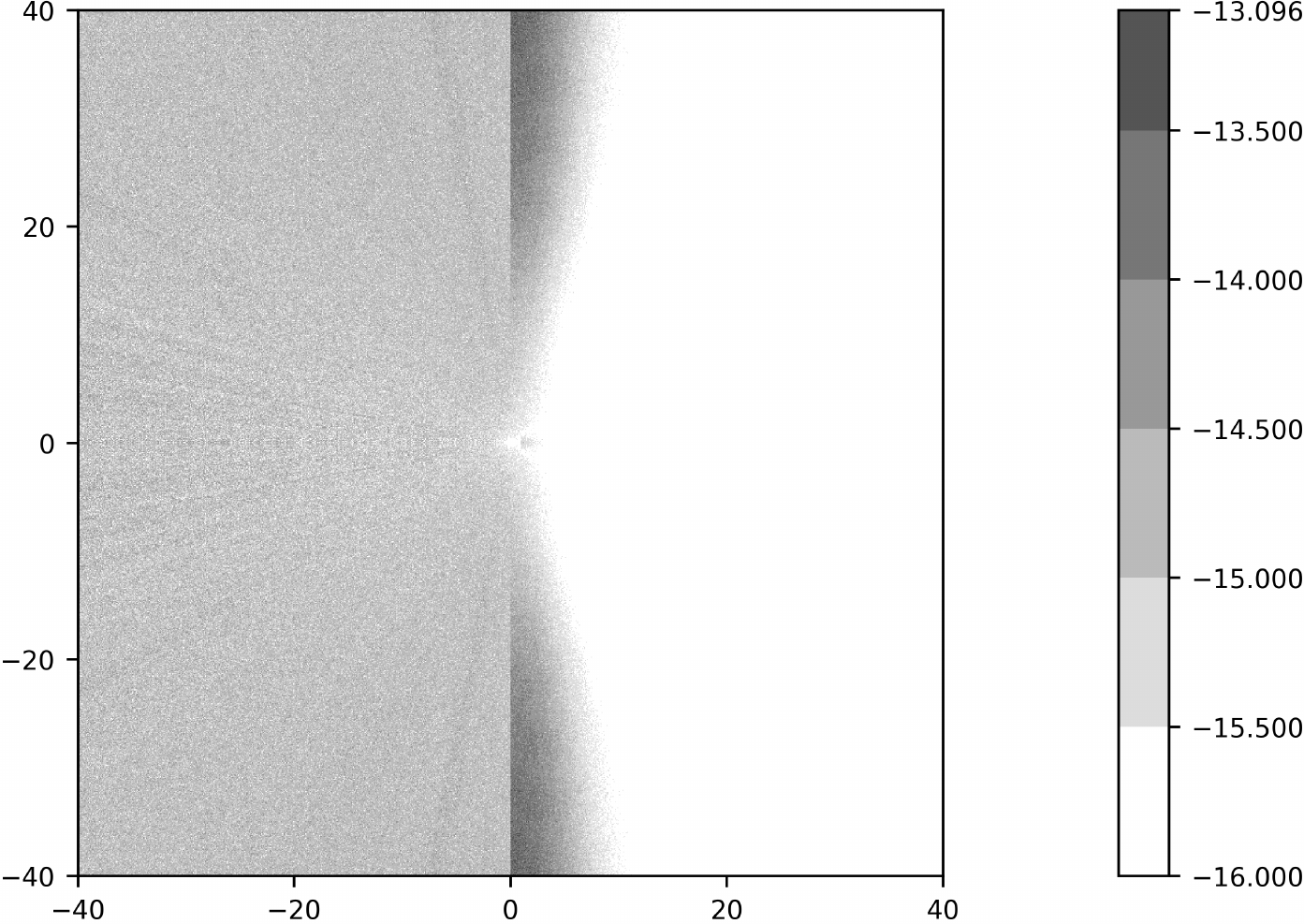}\kern25mm
 \includegraphics[height=46mm]{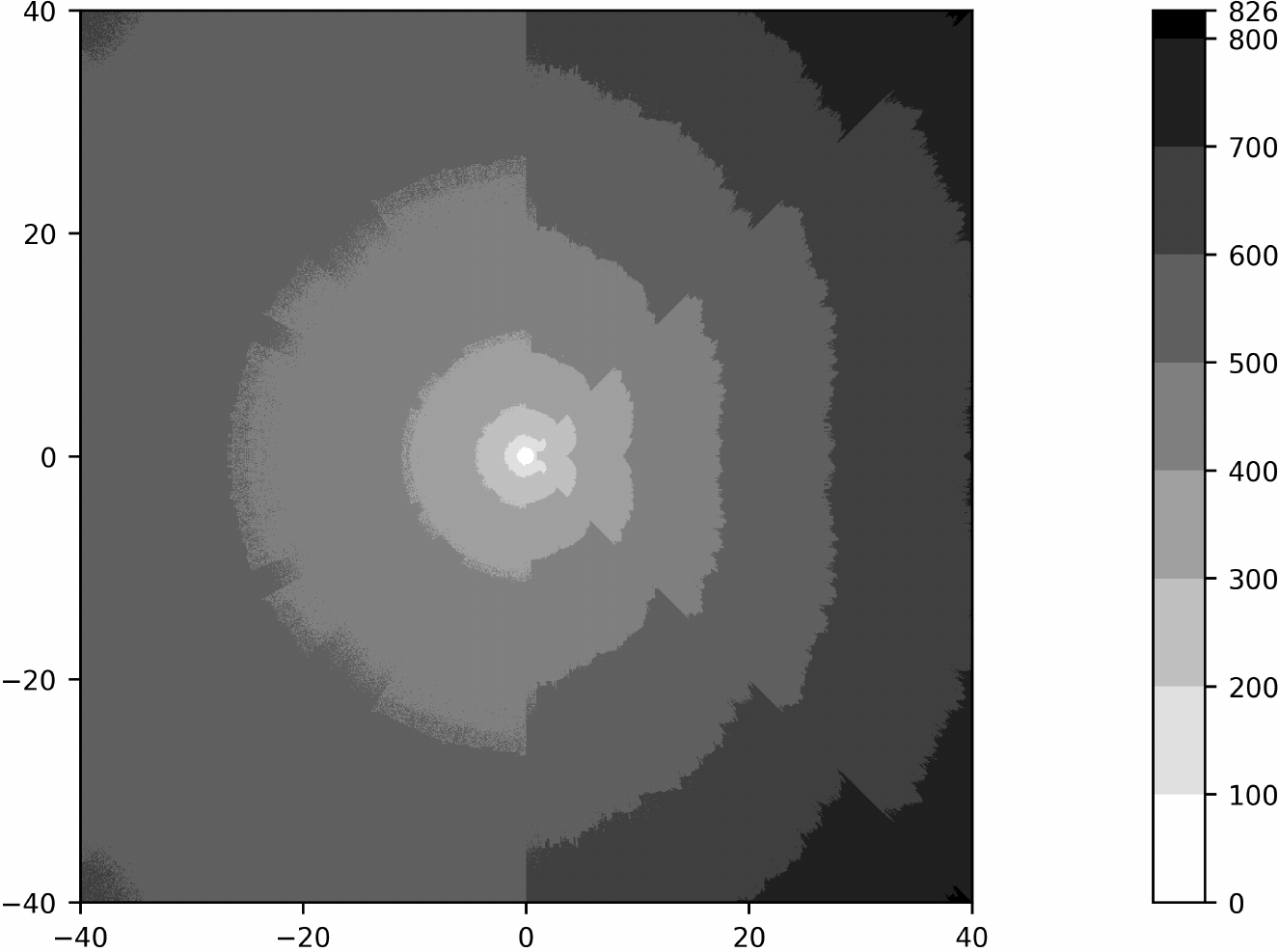}}
 \vspace{2mm}
 \caption{Values of $\max\bigl\{-16,\log_{10}\Lambda_8(z)\bigr\}$ (a) and $N_8(z)$ (b).}
 \label{fig:L8_0}
\centering\vspace{4mm}
 \SetLabels
 \L (-0.0425*0.96) (a)\\
 \L (0.521*0.96) (b)\\
 \L (0*0.86) $\Im z$\\
 \L (0.27*-0.06) $\Re z$\\
 \L (0.5675*0.86) $\Im z$\\
 \L (0.843*-0.06) $\Re z$\\
 \endSetLabels
 \leavevmode\AffixLabels{\mbox{\kern5mm}\includegraphics[height=46mm]{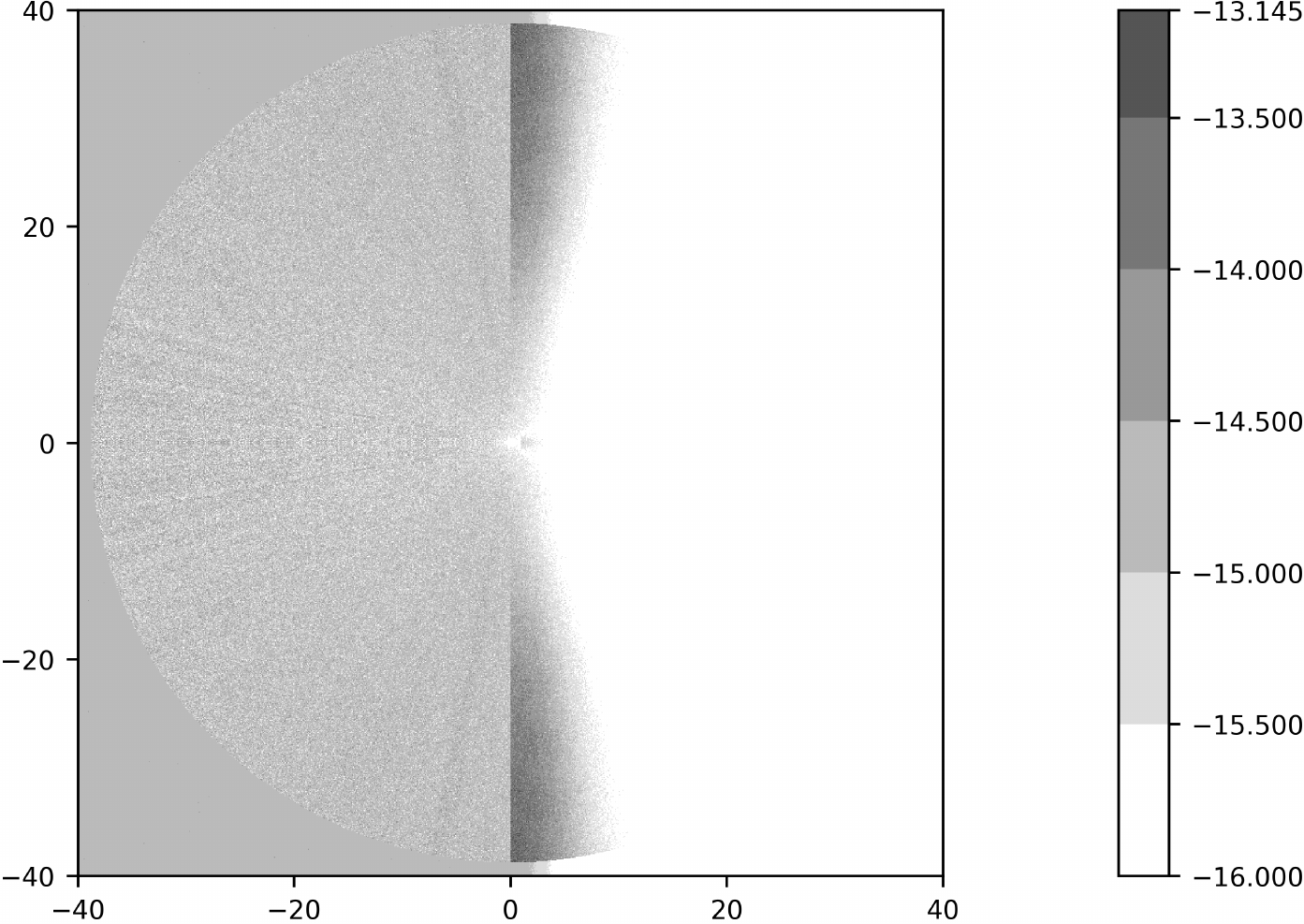}\kern25mm
 \includegraphics[height=46mm]{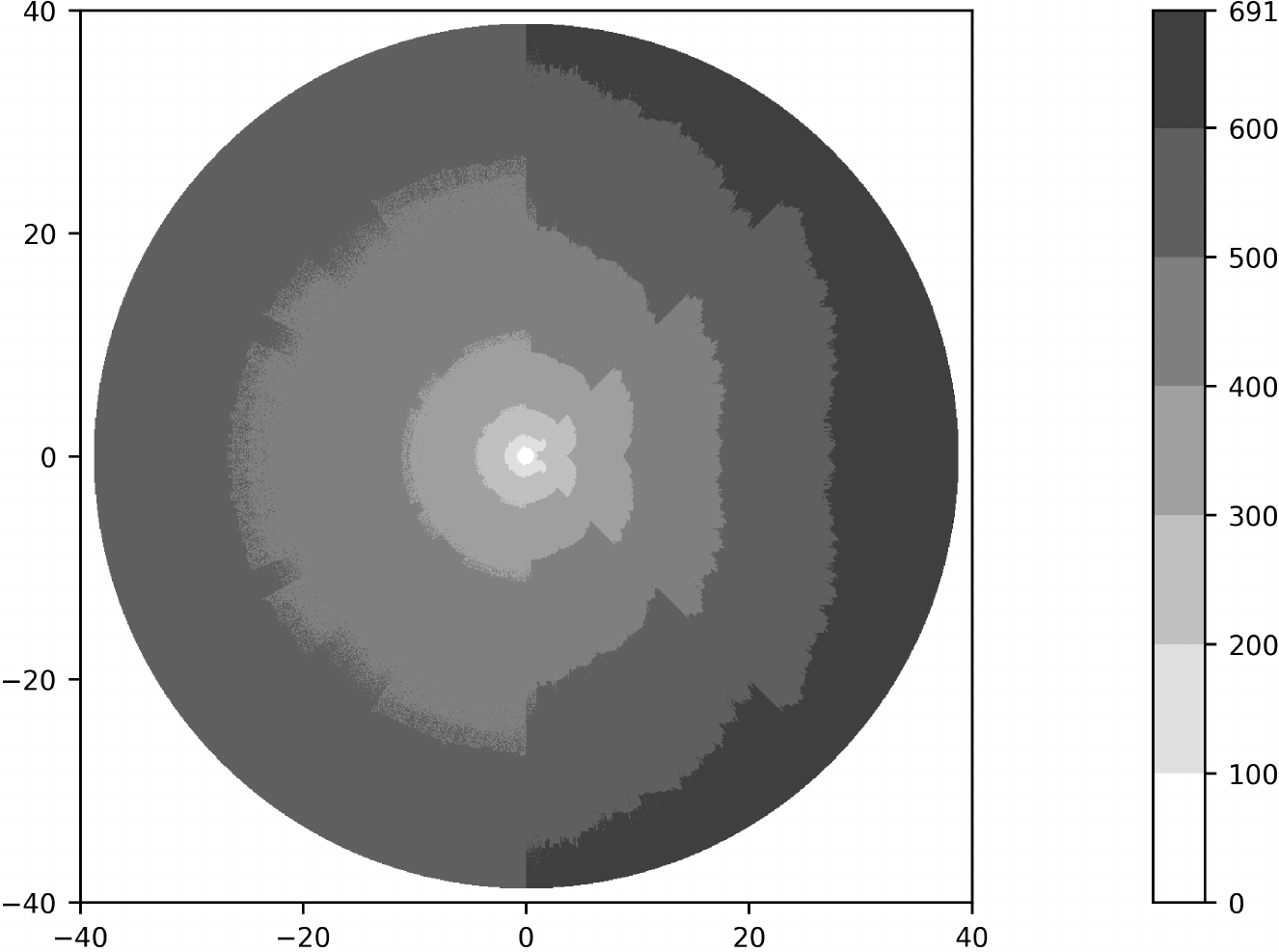}}
 \vspace{2mm}
 \caption{Values of $\max\bigl\{-16,\log_{10}\Lambda_8(z)\bigr\}$ (a) and $N_8(z)$ (b). Improved algorithms $\psHeunCsnearone$ and $\psHeunCsnearinf$ are used.}
 \label{fig:L8}
\end{figure}

\begin{figure}[t!]
\centering
 \SetLabels
 \L (-0.0425*0.96) (a)\\
 \L (0.521*0.96) (b)\\
 \L (0*0.86) $\Im z$\\
 \L (0.27*-0.06) $\Re z$\\
 \L (0.5675*0.86) $\Im z$\\
 \L (0.843*-0.06) $\Re z$\\
 \endSetLabels
 \leavevmode\AffixLabels{\mbox{\kern5mm}\includegraphics[height=46mm]{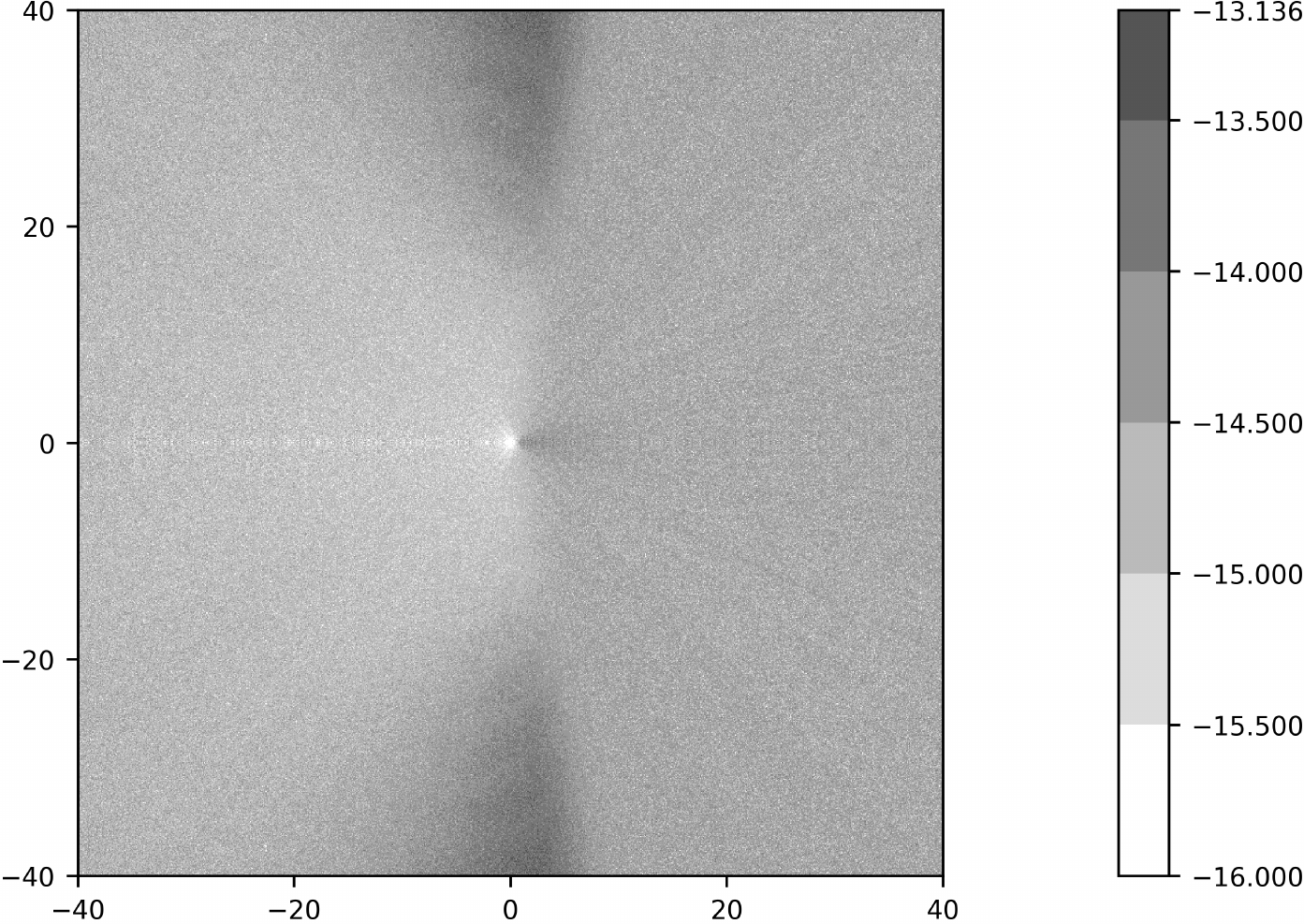}\kern25mm
 \includegraphics[height=46mm]{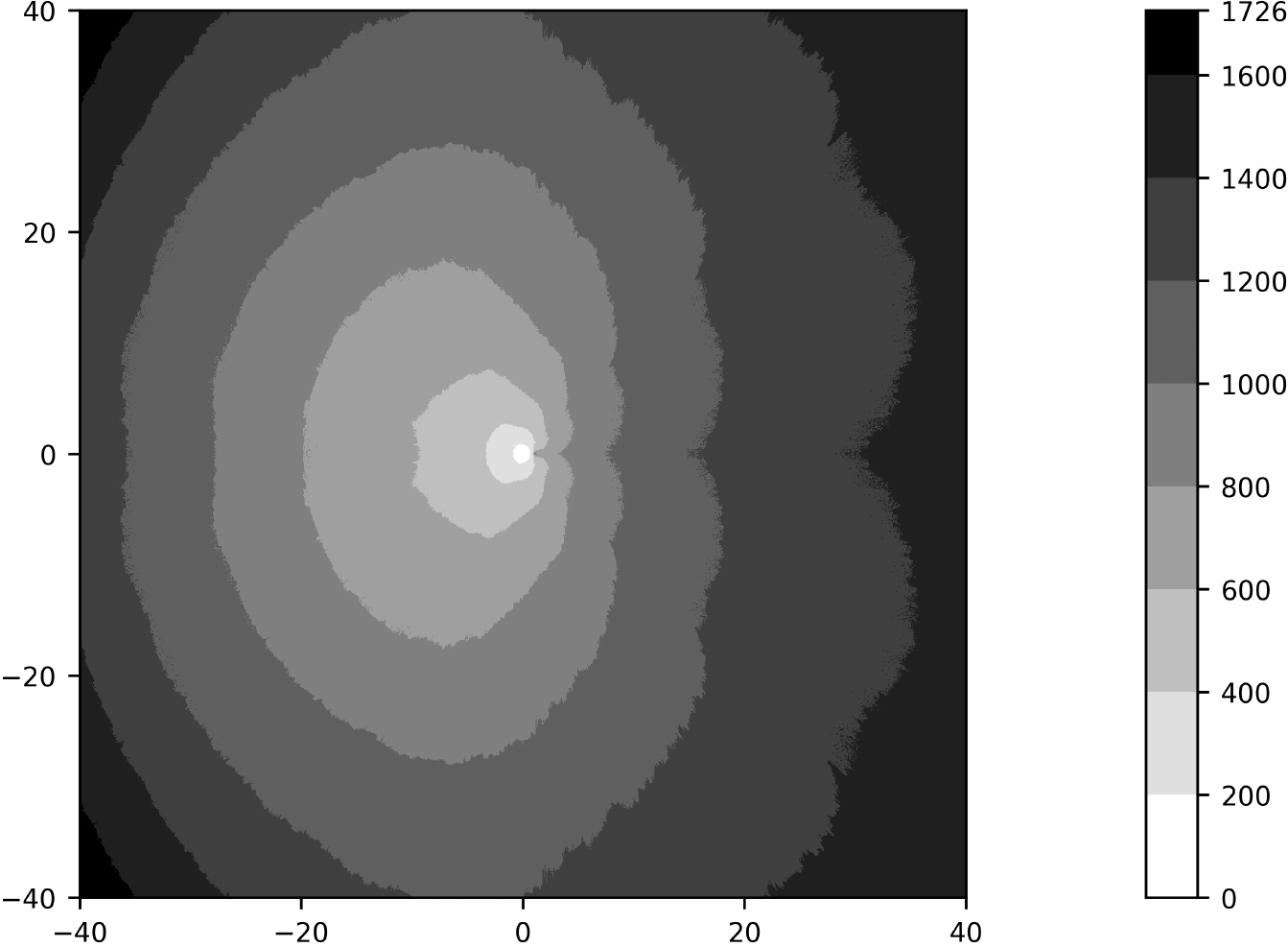}}
 \vspace{2mm}
 \caption{Values of $\max\bigl\{-16,\log_{10}\Lambda_9(z)\bigr\}$ (a) and $N_9(z)$ (b).}
 \label{fig:L9_0}
\end{figure}

Figures~\ref{fig:L1_0}a and \ref{fig:L2_0}a show in a semilogarithmic scale results of computations
of the relative error
\[
 \Lambda_n(z)=\frac{|\Delta_n(z)|}{1+|h_n(z)|}+\frac{|\Delta_n'(z)|}{1+|h_n'(z)|}
\]
for $n=1$, using the algorithm $\psHeunClfin$, and for $n=2$, using the algorithm $\psHeunCsfin$. 
In Fig.~\ref{fig:L1_0}b and \ref{fig:L2_0}b we present the values $N_1(z)$ and $N_2(z)$
which mean the total number of terms in power series used to compute
$\HeunCl(z)$ and $\HeunCs(z)$ in the expression of $\Lambda_1(z)$ and $\Lambda_2(z)$, respectively. 
These values can characterize the time of computation.

In these and other figures below, $\Lambda_n(z)$ and $N_n(z)$, $n=1,2,\ldots,9$, are computed on the grid $(\Re z,\Im z)\in\vec{L}(1000,[-40,40])\times\vec{L}(1000,[-40,40])$, where $\vec{L}(m,\chi)$ is the set of $m$ linearly spaced in the interval $\chi$ values (including interval's end-points). 

In Fig.~\ref{fig:L3_0} we present results of computations
of $\Lambda_3(z)$ and $N_3(z)$, using the algorithm $\psHeunClfin$. Fig.~\ref{fig:L4_0} shows the values 
of $\Lambda_4(z)$ and $N_4(z)$, computed by the algorithm $\psHeunCsfin$. We also applied the improved algorithm ($\psHeunClnearone$, $\psHeunCsnearone$); its accuracy is higher but pictures of computations using the improvement look very similar to Fig.~\ref{fig:L3_0} and \ref{fig:L4_0} and therefore omitted. 

Figures~\ref{fig:L5_0}a and \ref{fig:L6_0}a show some lost of accuracy and increase of $\Lambda_5(z)$, $\Lambda_6(z)$ above the real axis, but this presumably happens due to features of realizations of $\log$-functions in Octave, because in the computations $h_5(\overline{z})\neq\overline{h_5(z)}$ and $h_6(\overline{z})\neq\overline{h_6(z)}$.


In Figs.~\ref{fig:L7_0}--\ref{fig:L8} we compare the basic algorithm against the algorithm with improvements near points $z=1$ and $z=\infty$. Here  the improvement is more manifestative than for $\Lambda_3(z)$, $\Lambda_4(z)$. The values $N_7(z)$ and $N_8(z)$ shown in Figs.~\ref{fig:L7}, \ref{fig:L8} do not count operations needed for finding connection coefficients of local solutions at $z=0$ and $z=\infty$. When it is done, the algorithm uses  the coefficients saved in memory. Fig.~\ref{fig:L9_0} presents numerical results for $\Lambda_9(z)$ and $N_9(z)$. 

Finally, we note that Figs.~\ref{fig:L1_0}--\ref{fig:L9_0} do not show essential degradation of accuracy at increase of $|z|$, though, of course, computational load grows. The accuracy in the examples of computation is seemed to be fairly satisfactory for the used double-float arithmetics.

\end{document}